\documentclass[10pt]{article}
\usepackage{amsmath}
\usepackage{amssymb}
\usepackage{graphicx}
\usepackage{wrapfig}
\usepackage{float}
\usepackage{enumitem}
\usepackage{xcolor}
\usepackage{bbm} % Pour les ensembles de nombres
\usepackage{amsthm}
\usepackage{mathtools}
\usepackage{setspace}
\usepackage[left=2.5cm, right=2.5cm, top=2.5cm, bottom=1.5cm]{geometry}
\usepackage{bm} % bold in math mode
\usepackage{fancyhdr}
\usepackage[medium]{titlesec}
\usepackage{subfiles}
\usepackage{bigints} % grandes intégrales
\usepackage[autostyle=true]{csquotes} % guillemets
\usepackage[backend=biber, maxnames=5, minnames=1, sorting=nyt]{biblatex}
\newtheorem{theorem}{Theorem} % Crée l'environnement \begin{theorem}
\newtheorem{lemma}{Lemma} % Crée l'environnement \begin{lemma}
\newtheorem{proposition}{Proposition} % Crée l'environnement \begin{proposition}
\theoremstyle{remark}
\newtheorem{remark}{Remark}
\theoremstyle{definition}
\newtheorem*{definition}{Definition}

\newcounter{hypothese} % crée un nouveau compteur, appelé hypothese
\newcounter{hypspe} % crée un nouveau compteur, appelé hypspe
\renewcommand\thehypspe{\arabic{hypspe}} % le compteur hypspe sera affiché en lettres, mais attention, on le désigne toujours par un nombre, entre 1 et 26

\newenvironment{hypothese}[1][] % crée un nouvel environnement, appelé hypothese
{\refstepcounter{hypothese}\par \textbf{(H\thehypothese)#1} } % exécuté au début, avec le \begin{hypothese}
{\smallskip} % exécuté à la fin, avec le \end{hypothese}

{\refstepcounter{hypspe}\par \textbf{(H\thehypspe b)#1} } % exécuté au début, avec le \begin{hypspe}
{\smallskip} % exécuté à la fin, avec le \end{hypspe}

\usepackage{hyperref}

\hypersetup{
    colorlinks=true,
    linkcolor=black,    
    urlcolor=black,
    citecolor=black,
    }

\newcommand{\eps}{\varepsilon}
\newcommand{\R}{\mathbb{R}}
\newcommand{\N}{\mathbb{N}}
\newcommand{\ie}{\textit{i}.\textit{e}.}
\newcommand{\dd}{\mathrm{d}}
\newcommand{\ds}{\displaystyle}
\newcommand{\C}[1]{\mathcal{C}^{#1}}

\newcommand{\norm}[1]{\left\lVert #1\right\rVert}

\newcommand{\Cmin}{\underline{C}_{\min}}
\newcommand{\indic}[1]{\mathbbm{1}_{#1}}

\newcommand{\thetabar}{\bar{\theta}_d}

\title{Optimal Dirac controls for time-periodic bistable ODEs, application to population replacement}
\date{January 13th, 2026}
\author{Grégoire Nadin\footnote{CNRS, Institut Denis Poisson UMR 7013, Université d'Orléans, Collegium Sciences et Techniques, Bâtiment de mathématiques, Rue de Chartres, BP 6759, 45067 Orléans cedex 2, France (gregoire.nadin@univ-orleans.fr)}
\and David Nahmani\footnote{LAGA UMR 7539, Université Sorbonne Paris Nord, Institut Galilée, 99 avenue Jean-Baptiste Clément, 93430 Villetaneuse, France (nahmani@math.univ-paris13.fr)}
\and Nicolas Vauchelet\footnote{LAGA UMR 7539, Université Sorbonne Paris Nord, Institut Galilée, 99 avenue Jean-Baptiste Clément, 93430 Villetaneuse, France (vauchelet@math.univ-paris13.fr)}
}

\begin{document}

\maketitle

\begin{abstract}
This work addresses an optimal control problem on a dynamics governed by a nonlinear differential equation with a bistable time-periodic nonlinearity. This problem, relevant in population dynamics, models the strategy of replacing a population of $A$-type individuals by a population of $B$-type individuals in a time-varying environment, focusing on the evolution of the proportion of $B$-type individuals among the whole population.
The control term accounts for the instant release of $B$-type individuals. Our main goal, after noting some interesting properties on the differential equation, is to determine the optimal time at which this release should be operated to ensure population replacement while minimizing the release effort.
The results establish that the optimal release time appears to be the minimizer of a function involving the carrying capacity of the environment and the threshold periodic solution of the dynamics; they also describe the convergence of the whole optimal release strategy. An application to the biocontrol of mosquito populations using \textit{Wolbachia}-infected individuals illustrates the relevance of the theoretical results. \textit{Wolbachia} is a bacterium that helps preventing the transmission of some viruses from mosquitoes to humans, making the optimization of \textit{Wolbachia} propagation in a mosquito population a crucial issue.

\paragraph{Keywords: }Bistable differential equation, time-periodic differential equation, optimal control, optimization, population replacement, \textit{Wolbachia}

\paragraph{MSC: }34H05, 92D25, 49J15
\end{abstract}

\section{Introduction and main results} \label{main_results}

\maketitle

\subsection{Biological motivation}
The mathematical study of the present article was inspired by a population dynamics problem relevant in ecology. Assume that we want to replace a population of $A$-type individuals of any living species by a population of $B$-type individuals of the same species. This operation is called \textit{population replacement} and can be achieved releasing $B$-type individuals among the initial $A$-type population. For the sake of simplicity in the study, we only consider variations in time, and assume that this whole problem is observed at a given location, therefore the space-dependency is ignored in this paper.

Our objective is then to mathematically optimize this strategy, namely to determine when and how $B$-type individuals should be released, such that the replacement eventually occurs, and, above all, in the most economic way, considering that this process involves a financial cost associated with breeding and releases logistics.

First, we will stick to this general optimization problem, in order to come up with mathematical results with the largest possible scope. Then, further in the paper, we will detail an application to a model coming from a specific public health issue: a biocontrol strategy for mosquito populations.

The proportion of $B$-type individuals in the whole population at time $t$ is denoted $p(t)$ and follows an intrinsic dynamics governed by the nonlinear differential equation $p' = f(t,p)$.
We assume that the nonlinearity %$f(t,p)$
periodically depends on time, to take the seasonal variations of the habitat into account. Moreover, we consider nonlinearities of bistable type with respect to $p$. This assumption, which arises naturally in modeling population dynamics (see \cite{aronson-weinberger}, \cite{barton}, \cite{schraiber}), takes the competition between the two populations into consideration (see \cite{strugarek}) and induces a negative growth rate at low densities, $\ie$ a strong Allee effect (see \cite{wang}).

To know about controllability of bistable -- more general -- reaction-diffusion partial differential equations, involving space dependency, the reader can refer to \cite{affili} and \cite{pouchol}.

We add to this dynamics a control part of the form $u(t) w(t,p)$, where $u(t)$ denotes the released $B$-type population at time $t$, %scaled to $K(t)$, the environmental carrying capacity at time $t$,
and $w$ is a function that encompasses the time-dependency of the habitat and models the weight, depending on the value of $p$, of this control part on the overall equation.

Gathering the intrinsic dynamics and the control part, we study a control equation in which $f$ and $w$ satisfy some minimal hypotheses. The goal of this paper is to get information on the optimal released population in order to minimize a certain loss function.

Several mathematical works have been carried out on the optimization of such systems, and more specifically on models describing mosquito population dynamics. In all of them, the hosting capacity of the environment is assumed to be constant through time.

The work \cite{almeida19} studies the maximization of the proportion of $B$-type individuals at a fixed horizon of time, the total amount of released individuals through time being bounded. 
The same type of problems is spatially addressed in \cite{nadin}, \cite{mazari} and \cite{duprez}.
The article \cite{almeidachap} studies the minimization of the total amount of released individuals such that the critical equilibrium of the differential equation is reached at the final time, considering that the number of released individuals at a given time is bounded. In \cite{almeida24}, the authors focus, with an epidemiological model, on the transmission of diseases from mosquitoes to humans, and intend to minimize the total amount of infected humans.
Finally, the paper \cite{orozco-gonzales} considers a system of two equations, and deals with the minimization of the time of intervention and the total amount of released individuals, such that the system gets to the good basin of attraction.

The main mathematical result obtained in \cite{almeida19}, and later in \cite{almeidachap}, is that the optimal release function is bang-bang, and equal to its upper bound at the beginning of the time window, for as long as needed to reach the stated objective, and then equal to 0.

In the present work, we not only establish general results that are applicable to other periodic bistable phenomena, but we also consider that the hosting capacity is time-dependent, which adds a complexity to the problem.

To simplify the work, the profiles of the releases will be limited to scaled indicator functions, and asymptotically, to instant impulses. We will note further in the paper that this simplification has a biological relevance.

\subsection{Differential equation without control}

In this subsection, we investigate the ordinary differential equation
\begin{equation} \label{ODE}
    \forall t \in \R_+, \; p'(t) = f\bigl(t,p(t)\bigr),
\end{equation}
with $f : \R_+ \times [0,1] \to \R$ satisfying the following hypotheses:

\begin{hypothese} \label{hyp_regularite}
$f$ is $\C{2}$ in $\R_+ \times [0,1]$.
\end{hypothese}

\begin{hypothese} \label{hyp_periodicite}
$f$ is $T$-periodic with respect to the variable $t$.
\end{hypothese}

\begin{hypothese} \label{hyp_zeros}
For all $t \in \R_+, \; f(t,0) = f(t,1)=0$.
\end{hypothese}

\begin{hypothese} \label{hyp_vm_derivee_neg}
$\ds \frac{1}{T} \int_0^T \partial_p f\bigl(t,0\bigr) \dd t < 0 \text{ and } \ds \frac{1}{T} \int_0^T \partial_p f\bigl(t,1\bigr) \dd t < 0.$
\end{hypothese}

We already know that there exist trivial periodic solutions of the ODE \eqref{ODE}, namely $p\equiv 0$ and $p\equiv 1$.

Let us also introduce the mean value of $f$ over a period, in order to get an autonomous function:
$$\begin{array}{ccrcl}
f_m & : & [0,1]  & \to     & \R \\
  &   & p      & \mapsto & \ds \frac{1}{T} \int_0^T f(t,p) \dd t.
\end{array}$$
With the change of variables $s = t / T$, one gets $f_m(p) = \ds \int_0^1 \hat{f}(s,p) \dd s$, denoting $\hat{f}$ the 1-periodic function with respect to the first variable defined by $\hat{f}(s,p) := f(Ts,p)$. Note that the hypothesis (H\ref{hyp_vm_derivee_neg}) has an equivalent formulation, denoted \textbf{(H\ref{hyp_vm_derivee_neg}b)}, replacing $\partial_p f$ by $\partial_p \hat{f}$ and integrating between 0 and 1 instead of between 0 and $T$.

Let us now define the attractivity of a specific solution.

\begin{definition}
The $T$-periodic solution $p^T$ is said to be \emph{attractive} if there exists $\eta > 0$ such that for all $p$ solution of \eqref{ODE},
$$\left|p(0)-p^T(0)\right| < \eta \Longrightarrow \ds \lim_{t \to + \infty} \left|p(t) - p^T(t)\right| = 0.$$
\end{definition}

The positive number $\eta$ characterizes the basin of attraction of $p^T$.

B. Contri mentions in \cite{contri} (Propositions 2.1, 2.2 and 2.3) the \emph{principal eigenvalue} associated with the function $f$ and the $T$-periodic solution $p^T$, given by the expression
$$\lambda = - \frac{1}{T} \int_0^T \partial_p f\bigl(t,p^T(t)\bigr) \dd t.$$

\par Moreover, he establishes the following useful result: if $\lambda > 0$ (resp. $< 0$), then the $T$-periodic solution $p^T$ is attractive (resp. not attractive).
Given those tools, it is straightforward by (H\ref{hyp_vm_derivee_neg}) that the solutions $p \equiv 0$ and $p \equiv 1$ both are attractive.

Now let us state the following existence result, proved in \cite{contri} (Lemma 5.1).
\begin{proposition} \label{existence_periodic}
There exists at least one $T$-periodic solution $p^T:\R_+\to \; (0,1)$ of the ODE \eqref{ODE}.
\end{proposition}

Another important result is the existence of a maximal and a minimal $T$-periodic solutions.

\begin{proposition} \label{solutions_maximale_minimale}
Among the $T$-periodic solutions of the ODE \eqref{ODE} taking values in $(0,1)$, there exist a maximal solution $p^T_M$ and a minimal solution $p^T_m$.
\end{proposition}

The following proposition will be useful for the study of the optimal control problem.
It claims that any small perturbation outside the range of periodic solutions steers the dynamics to 0 or 1, the two attractive equilibria.

\begin{proposition} \label{basins_attraction}
Let $\bar{t} \in \R_+$.
Let $p$ be the solution of \eqref{ODE} that satisfies the condition $p(\bar{t}) = p_0$. Moreover, let us denote $p^T_M$ and $p^T_m$ respectively the maximal and the minimal $T$-periodic solutions strictly between 0 and 1.
\begin{enumerate}[label=(\roman*)]
\item If $p_0 > p^T_M(\bar{t})$, then $\ds \lim_{t\to +\infty} p(t) = 1$;
\item If $p_0 < p^T_m(\bar{t})$, then $\ds \lim_{t\to +\infty} p(t) = 0$.
\end{enumerate}
\end{proposition}

\begin{remark}
If $p^T_m(\bar{t}) \leqslant p_0 \leqslant p^T_M(\bar{t})$, then for all $t \geqslant \bar{t}$, one has $p^T_m(t) \leqslant p(t) \leqslant p^T_M(t)$, by Cauchy-Lipschitz theorem.
\end{remark}

\subsection{Optimal control problem} \label{main_results-OCP}

The control strategy we intend to study is population replacement, in other words replacing a population of $A$-type individuals by a population of $B$-type individuals. We will here aim the success of this strategy, with a mathematical point of view.

Let $p(t)$ denote the proportion of $B$-type individuals and $K(t)$ the environmental carrying capacity at time $t$. To analyze mathematically the success of population replacement strategy, let us investigate here a control equation, where one acts on the population, releasing a time-distributed amount of $B$-type individuals $u(t)$.

\paragraph{Control equation.}
We assume for now, the adequation with applications being detailed later, that this control equation is:
\begin{equation} \label{controlsyst_g_wo_K}
\begin{cases}    
    \ds p'(t) = f\bigl(t,p(t)\bigr) + \frac{u(t)}{K(t)} g\bigl(p(t)\bigr) \quad \text{in $(0,+\infty)$},\\
    p(0) = p_0,
\end{cases}
\end{equation}
with $f$ satisfying the hypotheses (H\ref{hyp_regularite}) to (H\ref{hyp_vm_derivee_neg}), $g:[0,1] \to \R$ satisfying the following hypothesis:

\begin{hypothese} \label{hyp_g}
$g$ is Lipschitz continuous, decreasing, $g(0) = 1$ and $g(1) = 0$,
\end{hypothese}

and $K : \R_+ \to \R$ satisfying the following hypothesis:
\begin{hypothese} \label{hyp_K}
$K$ is continuous and $T$-periodic in $\R_+$. Moreover, $\ds \min_{[0,T]} K > 0$.
\end{hypothese}

Indeed, since $K$ stands for the largest amount of population that the environment can host, given available resources, it seems relevant to assume that it depends on the season of the year. This is why $K$ is assumed to be a $T$-periodic function. 

Since the control $u$ can be discontinuous, the right-hand side is not continuous with respect to $t$, therefore the initial value problem \eqref{controlsyst_g_wo_K} needs to be read in the sense of Carathéodory. In this sense, the function $p$ is said to be a solution of \eqref{controlsyst_g_wo_K} if $p$ is absolutely continuous and satisfies for all $t \geqslant 0$,
$$p(t) = p_0 + \int_0^t \Bigl(f\bigl(s,p(s)\bigr) + \frac{u(s)}{K(s)} g\bigl(p(s)\bigr) \Bigr) \dd s.$$
Moreover, all the considered controls belong to $L^\infty(\R_+)$, therefore \eqref{controlsyst_g_wo_K} has a unique solution $p$.

\paragraph{Form of the control term.}
Given a total amount $C$ of available $B$-type individuals, one wants to release them all at once, in other words, one considers a release taking the form of the largest Dirac impulse allowed, operated at a certain time $t_0 > 0$. The study of Dirac controls mathematically allows us to understand the behavior of solutions and consists in a first step towards studying more diverse forms of control terms.

In order to approach Dirac controls, let us first study controls defined as follows: there exists $\eps > 0$ such that $\eps \ll 1$ and the control term writes, for all $t \geqslant 0$,
$$u^\eps(t) = \ds \frac{C}{\eps} \indic{[t_0-\eps,t_0]}(t).$$

As mentionned in \cite{almeida24}, \cite{almeida20} and \cite{duprez}, this form of control term is relevant since releases are often operated for a short time with respect to the time window considered.
Moreover, the numerical simulations that we will expose further in the paper will show that limiting the theoretical study to scaled indicator functions, and asymptotically to Dirac impulses, is a natural choice.

Let $p^\eps$ be the solution (in the sense of Carathéodory) of the control equation associated with the control $u^\eps$:
\begin{equation} \label{controlsyst_g_wo_K_eps}
\begin{cases}    
    \ds
    (p^\eps)'(t) = f\bigl(t,p^\eps(t)\bigr) + \frac{u^\eps(t)}{K(t)} g\bigl(p^\eps(t)\bigr) \quad \text{in $(0,+\infty)$},\\
    p^\eps(0) = p_0.
\end{cases}
\end{equation}

%{\color{orange}The reader notes that for $\eps$ fixed, $u^\eps$ is piecewise continuous and any solution $p^\eps$ is continuous, therefore Cauchy-Lipschitz theorem in each piece ensures the uniqueness.}
The optimal control problem that we want to address is the following, denoting $S := 1 / \eps$:
\begin{equation} \label{P_S} \tag{$\mathcal{P}_S$}
\boxed{
\begin{gathered}
\text{Find the time $t_0$ at which the release should be operated}\\
\text{in order to ensure that the effort $C$ is minimal and to}\\
\text{guarantee population replacement in problem \eqref{controlsyst_g_wo_K_eps}.}
\end{gathered}
}
\end{equation}
The reader notes that the objective function to minimize is the total amount $C$ of released individuals such that any greater amount ``guarantees population replacement", $\ie$ the solution of the control equation converges towards 1, whereas the works \cite{almeida19} and \cite{duprez}, on similar control strategies, consider as the objective function the difference between 1 and the state at the final time of observation, for a given amount of released individuals $C$.
In this paper, we tackle the issue of only steering the state to the good basin of attraction, to let the dynamics ensure by itself the convergence to 1. This question is closer to the ones addressed in \cite{almeidachap} and \cite{orozco-gonzales}.

\subsection{Asymptotic minimization problem}

In the following result, we study the asymptotic control equation when $\eps$ goes to 0.

\begin{proposition} \label{convergencep}
Within the above framework, let $p^\eps$ be a solution of the control equation \eqref{controlsyst_g_wo_K_eps} associated with the control $u^\eps$.
Then $(p^\eps)_{\eps>0}$ pointwise converges and its limit, denoted $p$, is the unique solution of the problem
\begin{equation} \label{Cauchy_IVP}
    \begin{cases}
    p'(t) = f\bigl(t,p(t)\bigr) \quad \text{for all } t \in (0,t_0) \cup (t_0,+\infty),\\
    G\bigl(p(t_0^+)\bigr) - G\bigl(p(t_0^{-})\bigr)  = \ds \frac{C}{K(t_0)},\\
    p(0) = p_0,
    \end{cases}
\end{equation}
with $G$ an antiderivative of $1 / g$.
\end{proposition}

\begin{remark}
A proof of this result has been carried out in \cite{nedeljkov}, Proposition 2.1, and adapted in \cite{almeida24}, Proposition 1, for a population replacement problem. Our proof, in Subsection \ref{proof_proposition_cv_p}, uses a different method, first establishing the pointwise convergence relying on Helly's selection theorem, and then justifying the jump with a rescaling, whereas the one in both cited papers compares the solution to the solution of a simpler equation and states that their jumps at $t_0$ are the same.
\end{remark}

\paragraph{Heuristics on the asymptotic minimization problem.}
One now wants to solve \eqref{P_S} for the asymptotic control equation, $\ie$ when $\eps$ vanishes to 0. Assume that the system starts at 0, $\ie$ that $p_0 = 0$, which corresponds to a natural situation prior to population replacement, since at time $t=0$, there are only $A$-type individuals in the environment. Therefore, since the system between 0 and $t_0^-$ is simply the system without control, uniqueness in Cauchy-Lipschitz theorem yields $p \equiv 0$ in $(0,t_0)$.

Choosing for $G$ the antiderivative that cancels at 0, the asymptotic problem \eqref{Cauchy_IVP} reads:
\begin{equation} \label{Cauchy_IVP_0}
    \begin{cases}
    p \equiv 0 \text{ in } (0,t_0),\\[5pt]
    G\bigl(p(t_0^+)\bigr)  = \ds \frac{C}{K(t_0)},\\[10pt]
    p'(t) = f\bigl(t,p(t)\bigr) \quad \text{for all } t \in (t_0,+\infty).    
    \end{cases}
\end{equation}

According to Proposition \ref{basins_attraction}, $\ds \lim_{t\to +\infty}p(t) = 1$ in problem \eqref{Cauchy_IVP_0} is equivalent to: 
$$
\begin{array}{rcll}
p(t_0^+) > p^T_M(t_0) & \Longleftrightarrow & G\bigl(p(t_0)\bigr) > G\bigl(p^T_M(t_0)\bigr) & \text{since $G$ is increasing, by (H\ref{hyp_g}), therefore invertible,}\\[10pt]
    
    & \Longleftrightarrow & C > K(t_0) G\bigl(p^T_M(t_0)\bigr) &\text{since $K$ is positive,}\\[10pt]

    & \Longleftrightarrow & C > \underline{C}(t_0)
                        & \text{setting } \underline{C}(t_0) := K(t_0)  G\bigl(p^T_M(t_0)\bigr).\\ 
                    
\end{array}
$$

The condition is then to release at least $\underline{C}(t_0)$ $B$-type individuals at the time $t_0$. The optimal control problem \eqref{P_S} for the asymptotic equation seems therefore to be finding the time $t_0$ for which $\underline{C}(t_0)$ is minimal.

\paragraph{Rigorous result.}
Denote $S := 1 / \eps$. We intend to solve \eqref{P_S} for the asymptotic control equation, $\ie$ when $S$ goes to $+\infty$.

Let us introduce the following minimization problem:
\begin{equation} \label{P_eq} \tag{$\mathcal{P}_{\text{eq}}$}
\boxed{
\ds \inf_{t \in [0,T]} \underline{C}(t) \quad \text{with } \underline{C}(t) := K(t) G\bigl(p^T_M(t)\bigr).
}
\end{equation}

Note that the function $\underline{C}$ is continuous, thus the minimization problem \eqref{P_eq} has a solution in $[0,T]$.
The following theorem rigorously details the asymptotic minimization problem.

\begin{theorem} \label{OCP_asymptotic}
For $S > 0$ and $t_0^S \geqslant 1/S$, let us introduce the following control equation, for a given $C>0$:
\begin{equation} \label{S_SCt0S} \tag{$\mathcal{S}_{S,C,t_0^S}$}
	\begin{cases}
    p_S'(t) = f\bigl(t,p_S(t)\bigr) + \ds \frac{C S}{K(t)} \indic{[t_0^S-\frac{1}{S},t_0^S]}(t) \; g\bigl(p_S(t)\bigr) \quad\text{for almost every } t\in (0,+\infty),\\
    p_S(0) = 0,
    \end{cases}
\end{equation}
with $f$ satisfying (H\ref{hyp_regularite}) to (H\ref{hyp_vm_derivee_neg}), $g$ satisfying (H\ref{hyp_g}) and $K$ satisfying (H\ref{hyp_K}).

Set $\mathcal{I}_S := \Bigl\{C>0:\exists t_0^S \geqslant 1/S \;\text{ s. t. } p_S\text{ solution of \eqref{S_SCt0S} satisfies } 
\ds \lim_{t \to + \infty} p_S(t) = 1 \Bigr\}$.

First, for $S$ large enough, the set $\mathcal{I}_S$ has an infimum, denoted $C_S^*$.

Moreover, denoting $G$ the antiderivative of $1 / g$ vanishing at 0 and $$C^* := \ds \inf_{t \in [0,T]} \Bigl[K(t) G\bigl(p^T_M(t)\bigr)\Bigr] = K(t_0^*) G\bigl(p^T_M(t_0^*)\bigr),$$ one gets that $C_S^*$ converges towards $C^*$, when $S$ goes to $+\infty$.

Finally, the sequence of minimizers $\left(t_0^{S*}\right)_S$ of the set $\mathcal{I}_S$ converges, up to an extraction, to a global minimizer of $K(\cdot)G\bigl(p^T_M(\cdot)\bigr)$.
\end{theorem}

\begin{remark}
If the function $K(\cdot)G\bigl(p^T_M(\cdot)\bigr)$ has a unique global minimizer, then the solution(s) of \eqref{P_S} converge(s) to the unique solution of \eqref{P_eq} when $S$ goes to $+\infty$.
\end{remark}

\begin{remark}
The convergences of the sequence of infima and of the sequence of minimizers can be proved with $\Gamma$-convergence (see \cite{almeida24}, \cite{almeida19}). But here, since the asymptotic optimal control problem is actually a minimization problem in a real-valued interval, we use elementary and direct strategies to prove them (see Subsection \ref{proof_theorem}).
\end{remark}

\subsection{Uniqueness of the periodic solution}

In this subsection, we detail three situations in which the $T$-periodic solution strictly between 0 and 1 is unique. In those cases, one has $p^T_M = p^T_m = p^T$.

\subsubsection{Small periods}

Firstly, we consider the regime of small periods. In order to understand the limit as the period goes to 0, we introduce the 1-periodic function $\hat{f}$. The following result is established in \cite{contri} (Theorem 1.3, Subsection 5.1).
\begin{proposition} \label{uniqueness}
Let $\hat{f}$ be a 1-periodic function satisfying the hypotheses (H\ref{hyp_regularite}), (H\ref{hyp_zeros}), (H\ref{hyp_vm_derivee_neg}b) and, denoting $f_m(p) = \ds \int_0^1 \hat{f}(s,p) \dd s$,  the following hypothesis:
\begin{hypothese} \label{hyp_bistable_fm}
 There exists a unique $\bar{\theta} \in (0,1)$ such that $f_m(0) = f_m(\bar{\theta}) = f_m(1) = 0$, and one has $f_m < 0$ in $(0,\bar{\theta})$ and $f_m > 0$ in $(\bar{\theta},1)$, $\ie$ $f_m$ is bistable. Moreover, one has $f_m'(\bar{\theta})>0$.
\end{hypothese}

Let $f^T$ be the $T$-periodic function defined by $f^T(t,p):=\hat{f}(\frac{t}{T},p)$.

Then, there exists $T_{\hat{f}} > 0$ (only depending on $\hat{f}$) such that for all $T \in (0,T_{\hat{f}})$, the $T$-periodic solution $p^T$ of the ODE $p' = f^T(t,p)$ such that $0 < p^T(\cdot) < 1$ is unique and not attractive.
\end{proposition}

\subsubsection{Perturbed nonlinearity}

Secondly, we consider the situation of a perturbed nonlinearity. B. Contri establishes the following proposition in \cite{contri} (Theorem 1.5 and Subsection 6.1), relying on the mathematical framework of \cite{alikakos-bates-chen}.
\begin{proposition} \label{perturbed}
Let $f^0$ be a function that satisfies the hypotheses (H\ref{hyp_regularite}) to (H\ref{hyp_vm_derivee_neg}) and let $(f^\eps)_{\eps > 0}$ be a family of functions that all satisfy the same hypotheses, and such that there exists a bounded function $\omega :\R_+^* \to  \R$ such that $\omega(\eps) \xrightarrow[\eps \to 0]{}0$ and for all $\eps > 0$, one has
\begin{equation} \label{condition_f_der_eps}
\forall (t,p) \in [0,T] \times [0,1],\quad \big|\partial_p f^0(t,p) - \partial_p f^\eps(t,p) \big| \leqslant \omega(\eps).
\end{equation}
Assume that the ODE $p' = f^0(t,p)$ has a unique periodic solution $p^T$ strictly between 0 and 1 and that the eigenvalue associated with $f^0$ and $p^T$ is negative, $\ie \; \frac{1}{T} \int_0^T \partial_p f^0\bigl(t,p^T(t)\bigr) \dd t > 0$.
Then there exists $\eps_0 > 0$ such that for all $\eps \in (0,\eps_0)$, the perturbed ODE $p' = f^\eps(t,p)$ also has a unique periodic solution between 0 and 1 and this unique periodic solution is not attractive.
\end{proposition}

\subsubsection{Nonlinearity with separated variables}

Thirdly, we focus on a specific situation where the time-dependent part and the $p$-dependent part are separated in the expression of the nonlinearity, which gives the following result, that we establish, and whose proof is detailed in Subsection \ref{periodic_solutions}.

\begin{proposition} \label{separation_var_unicite}
Let $f : \R_+ \times [0,1] \to \R$ be defined by $f(t,p) := m(t) \tilde{f}(p)$, with $m$ and $\tilde{f}$ satisfying the following hypotheses:

\begin{hypothese} \label{hyp_periodicite_bis}
$m : \R_+ \to \R $ is $\C{1}$, $T$-periodic and $\int_0^T m(t) \dd t > 0$.
\end{hypothese}

\begin{hypothese} \label{hyp_regularite_bis}
$\tilde{f} : [0,1] \to \R$ is $\C{2}$.
\end{hypothese}

\begin{hypothese} \label{hyp_vm_derivee_neg_bis}
$\tilde{f}'(0) < 0 \text{ and } \tilde{f}'(1) < 0$.
\end{hypothese}

\begin{hypothese} \label{hyp_bistable_bis}
There exists a unique $\bar{\theta} \in (0,1)$ such that $\tilde{f}(0) = \tilde{f}(\bar{\theta}) = \tilde{f}(1) = 0$, $\tilde{f} < 0$ in $(0,\bar{\theta})$ and  $\tilde{f} > 0$ in $(\bar{\theta},1)$, that is, $\tilde{f}$ is bistable. Finally, one has $\tilde{f}'(\bar{\theta}) > 0$.
\end{hypothese}

Then the ODE \eqref{ODE} has a unique $T$-periodic solution $p^T$ strictly between 0 and 1, which is constantly equal to $\bar{\theta}$. Moreover, this solution is not attractive.
\end{proposition}

\begin{remark}
The reader can notice that such a function $f$ satisfies hypotheses (H\ref{hyp_regularite}) to (H\ref{hyp_vm_derivee_neg}) and (H\ref{hyp_bistable_fm}).
\end{remark}

\subsection{Optimal problem when the nonlinearity has separated variables}

In this subsection, we place ourselves in the framework of a nonlinearity with separated variables satisfying Proposition \ref{separation_var_unicite}. This simplified framework models the independence of effects of time and of the population proportion $p$ on the nonlinearity $f$.
By Proposition \ref{separation_var_unicite}, we already know that there exists a unique periodic solution strictly between 0 and 1 and that it is constantly equal to $\bar{\theta}$.
This favorable framework allows us to rewrite the optimal control problem \eqref{P_eq} as follows:
\begin{equation*}
\ds \inf_{t \in [0,T]} \Bigl[K(t) G(\bar{\theta}) \Bigr],
\end{equation*}
with $G$ the antiderivative of $1 / g$ vanishing at 0.

Let us consider first a single release at the time $t_0$. This optimization problem is equivalent to minimizing $K$, since $G$ is positive in $(0,1)$. Let us denote $\Cmin := \ds G(\bar{\theta})\min_{t \in [0,T]} K(t)$.
The best time $t_0^*$ to release the whole amount of $B$-type individuals is the time when the carrying capacity $K$ is minimal, and $\Cmin = \underline{C}(t_0^*)$, with $t_0^*$ a minimizer of $K$ in $[0,T]$.

The following proposition establishes a result on operating two releases.
Adapting the proof of Proposition~\ref{convergencep}, one can establish that, in the case of two releases, a solution $p^\eps$ of the control equation
\begin{equation} \label{controlsyst_g_wo_K_eps_2lachers}
\begin{cases}    
    \ds
    (p^\eps)'(t) = f\bigl(t,p^\eps(t)\bigr) + \frac{u_0^\eps(t) + u_1^\eps(t)}{K(t)} g\bigl(p^\eps(t)\bigr) \quad \text{for almost every $t \in (0,+\infty)$},\\
    p^\eps(0) = 0,
\end{cases}
\end{equation} 
with $u_0^\eps(t) = \ds \frac{C_0}{\eps} \indic{[t_0-\eps,t_0]}(t)$ and $u_1^\eps(t) = \ds \frac{C_1}{\eps} \indic{[t_1-\eps,t_1]}(t)$,
pointwise converges towards the solution $p$ of
\begin{equation} \label{Cauchy_IVP_2lachers}
    \begin{cases}
    p'(t) = f\bigl(t,p(t)\bigr) \quad \text{for all } t \in (0,t_0) \cup (t_0,t_1) \cup (t_1,+\infty),\\
    G\bigl(p(t_0^+)\bigr) - G\bigl(p(t_0^{-})\bigr)  = \ds \frac{C_0}{K(t_0)},\\[8pt]
    G\bigl(p(t_1^+)\bigr) - G\bigl(p(t_1^{-})\bigr)  = \ds \frac{C_1}{K(t_1)},\\
    p(0) = 0,
    \end{cases}
\end{equation}
with $G$ an antiderivative of $1 / g$.

\begin{proposition} \label{2releases}
Assume that the nonlinearity writes $f(t,p) = m(t) \tilde{f}(p)$ and satisfies Proposition \ref{separation_var_unicite}. Moreover, assume that for all $t > 0, \; m(t) > 0$.
Denote $t_0^*$ a minimizer of $K$ in $[0,T]$, $G$ the antiderivative of $1 / g$ vanishing at 0 and $\Cmin := \ds G(\bar{\theta}) K(t_0^*)$. Denote $t_0, t_1$ two fixed times such that $t_1 > t_0$ and consider $p$, the solution of the system \eqref{Cauchy_IVP_2lachers}, with $C_0, C_1 > 0$ such that $C_0+C_1 = \Cmin$.
One has $\ds \lim_{t \to + \infty} p(t) = 0$.
\end{proposition}

\begin{remark}
In other words, at the time $t_0$, we operate a first instant release of $C_0$ individuals and at the time $t_1$, we operate a second instant release of $C_1$
individuals.
This strategy of two releases will not lead to the basin of attraction of 1 and is consequently less efficient than operating a single release of $\Cmin$ individuals at the time $t_0^*.$
We could ask ourselves if this interesting practical result holds in the general framework, without the separated variables hypothesis. This open problem might be addressed in a future work.
\end{remark}

\subsection{Outline of the article}
The main results on the differential equation without control, such as Propositions \ref{solutions_maximale_minimale}, \ref{basins_attraction}, and \ref{separation_var_unicite}, will be established in Section \ref{proofs_1}. The main results on the optimal control equation, such as Propositions \ref{convergencep} and \ref{2releases} and Theorem \ref{OCP_asymptotic}, will be proved in Section \ref{proofs_2}\footnote{Propositions \ref{existence_periodic}, \ref{uniqueness} and \ref{perturbed} are proved in \cite{contri}.}. Section \ref{application_wolbachia} is dedicated to applying those mathematical results to a specific control strategy for mosquito populations, that has already been successfully implemented in some regions: population replacement based on the bacterium \textit{Wolbachia}.

% Ne pas insérer de bibliographie à la fin de chaque sous-fichier

\section{Differential equation without control} \label{proofs_1}

\maketitle

In this section, we investigate the ordinary differential equation \eqref{ODE} with $f$ satisfying the hypotheses (H\ref{hyp_regularite}) to (H\ref{hyp_vm_derivee_neg}).

First, hypotheses (H\ref{hyp_regularite}) and (H\ref{hyp_periodicite}) imply that $\partial_p f$ is continuous, and it is $T$-periodic with respect to $t$, thus bounded in $\R_+ \times [0,1]$, therefore $f$ is Lipschitz continuous with respect to $p$, which allows us to use Cauchy-Lipschitz theorem.
In particular, it gives that, knowing that $p \equiv 0$ and $p \equiv 1$ are solutions of \eqref{ODE}, if $p$ denotes a solution of \eqref{ODE} such that $p(0) \in [0,1]$, then for all $t \in \R_+, \; p(t) \in [0,1]$.

\subsection{Maximal and minimal periodic solutions, proof of Proposition \ref{solutions_maximale_minimale}}

\begin{definition}
For any $p_0 \in [0,1]$, let $p(\, \cdot \, ; p_0)$ be the solution of the Cauchy problem
\begin{equation*}
\begin{cases}    
    \ds p'(t) = f\bigl(t,p(t)\bigr) \quad \text{for all $t \geqslant 0$},\\
    p(0) = p_0.
\end{cases}
\end{equation*}
The \emph{Poincaré map} associated with $f$ is the function $P : [0,1] \to [0,1]$ defined by $P(p_0) := p(T ; p_0)$.
\end{definition}

\begin{proof}[Proof of Proposition \ref{solutions_maximale_minimale}]
Let $P$ be the Poincaré map associated with $f$, and $\Phi := P - \mathrm{Id}$; denote $Z := \{x \in (0,1) : \Phi(x) = 0\}$.
Proposition 2.4 in \cite{contri} states that, for $\alpha$ a fixed point of $P$, $P'(\alpha) = e^{-\lambda T}$, with $\lambda = - \frac{1}{T} \int_0^T \partial_p f\bigl(t,p(t;\alpha) \bigr) \dd t$.

Since $\Phi(0) = \Phi(1) = 0$, $\Phi'(0) < 0$, $\Phi'(1) < 0$, by \cite{contri} (Propositions 2.2, 2.3 and 2.4), and $\Phi'$ is continuous in $[0,1]$, then there exist $b$ and $c$ in $(0,1)$ such that $\Phi < 0$ in $(0,b)$ and $\Phi > 0$ in $(c,1)$.
Therefore, one can rewrite the set $Z$ as $Z = \{x \in [b,c] : \Phi(x) = 0\}$, and, since $\Phi : [b,c] \to \R$ is continuous, $Z$ is a closed subset of $[b,c]$, hence of $\R$.
Moreover, as a non-empty and bounded from above set, $Z$ has a supremum, which is a maximum, denoted $x_M$.

It is straightforward that the zeros of $\Phi$ are the initial values of $T$-periodic solutions. Therefore, the solution $p^T_M$ starting at $x_M$ is $T$-periodic, and by uniqueness in Cauchy-Lipschitz theorem, it stays strictly between 0 and 1 and remains greater than any other $T$-periodic solution, $\ie$ it is maximal.

With the same argument, there exists $x_m \in [b,c]$ such that the solution $p^T_m$ starting at $x_m$ is minimal among the $T$-periodic solutions strictly between 0 and 1.
\end{proof}

\subsection{Basins of attraction of the equilibria, proof of Proposition \ref{basins_attraction}}

\begin{proof}[Proof of Proposition \ref{basins_attraction}]
Let us denote $\overline{P}$ the Poincaré map shifted at time $\bar{t}$, $\ie$:
$$
\begin{array}{ccrcl}
\overline{P} & : & [0,1] & \to     & [0,1]\\
  &   &  p_0  & \mapsto & p(\bar{t}+T ; p_0),
\end{array}
$$
with $p(t ; p_0)$ denoting, for all $t \in [\bar{t},\bar{t}+T]$, the evaluation at $t$ of the solution $p$ satisfying the condition $p(\bar{t}) = p_0$.
The properties established in \cite{contri}, Section 2 hold for this shifted Poincaré map.
Denoting $\overline{\Phi} := \overline{P} - \mathrm{Id}$, one knows that $\overline{\Phi}(0) = \overline{\Phi}(1) = 0$, $\overline{\Phi}'(0) < 0$ and $\overline{\Phi}'(1) < 0$, by \cite{contri}, Propositions 2.2, 2.3 and 2.4.

Since $p^T_M$ and $p^T_m$ are respectively the maximal and the minimal $T$-periodic solutions taking values in $(0,1)$, $p^T_M(\bar{t})$ and $p^T_m(\bar{t})$ are respectively the maximal and the minimal zeros of $\overline{\Phi}$ in $(0,1)$. It yields $\overline{\Phi} < 0$ in $\bigl(0,p^T_m(\bar{t})\bigr)$ and $\overline{\Phi} > 0$ in $\bigl(p^T_M(\bar{t}), 1\bigr)$, by continuity of $\overline{\Phi}$ on $[0,1]$.

\begin{enumerate}[label=(\roman*)]
\item \label{1ercas} Let us assume that $p_0 > p^T_M(\bar{t})$. One gets $p(\bar{t}+T) - p(\bar{t}) = \overline{P}(p_0) - p_0 = \overline{\Phi}(p_0) > 0$, which gives $p(\bar{t}+T) > p(\bar{t})$.

The function $p(\, \cdot + T)$ satisfies the same equation as $p$. We can apply the same inequality to the function $p(\, \cdot + T)$, and we obtain similarly, since $p(\bar{t}+T) > p(\bar{t}) = p_0 > p^T_M(\bar{t})$,
$$p(\bar{t} + 2T) > p(\bar{t} + T).$$

For all $n \in \N^*$, we can show by induction $p(\bar{t} + nT) > p\bigl(\bar{t} + (n-1)T\bigr)$.
The sequence $\Bigl(p\bigl(\bar{t} + nT\bigr)\Bigr)_{n\in\N}$ is increasing and bounded, thus it converges to a certain limit $\overline{p^\infty}$. Moreover, one knows that for all $n \in \N, \; P\bigl(p(\bar{t}+nT)\bigr) = p\bigl(\bar{t}+(n+1)T\bigr)$, thus $\overline{p^\infty}$ is a fixed point of $P$.
By uniqueness in Cauchy-Lipschitz theorem, the solution $p$ cannot cross $p^T_M$, hence, $\overline{p^\infty} = 1$.

One understands that $\bar{t}$ can be chosen arbitrarily and the conclusion $p(\bar{t}+nT) \xrightarrow[n \to \infty]{} 1$ still holds.

Now, let $\eps > 0$ and $m \in \N^*$ such that $m\eps > T$, and let $l=T/m \in \R_+^*$.
For all $i \in \{0,...,m\}$, set $t_i = il$. The fact that $p(t_i+nT)$ tends to 1 implies that there exists $N_i \in \N$ such that for all $n \geqslant N_i,\; p(t_i+nT) > 1 - \eps$.

Let $N = \ds \max_{i\in\{0,...,m\}} N_i$ and $t > NT$. Setting $n = \lfloor \frac{t}{T}\rfloor$, one gets $n \geqslant N$ and $t-nT \in [0,T)$.
Thus, with $i = \lfloor \frac{t-nT}{l} \rfloor$, one has $t-nT \in [t_i,t_{i+1})$, $\ie$ $0 \leqslant t-(t_i+nT) < l < \eps$.

Finally, as one knows that $p'$ is bounded, the mean value theorem gives, for all $t > NT$,
$$\bigl|p(t) - p(t_i + nT)\bigr| < \norm{p'}_{\infty} \eps,$$ which implies
$$1 > p(t) > p(t_i + nT) - \norm{p'}_{\infty} \eps > 1 - (\norm{p'}_{\infty}+1) \eps.$$
Therefore $\ds \lim_{t\to+\infty}p(t) = 1$.

\item \label{2ecas} Let us now assume that $p_0 < p^T_m(\bar{t})$. One gets $p(\bar{t}+T) - p(\bar{t}) = \overline{P}(p_0) - p_0 = \overline{\Phi}(p_0) < 0$, which gives $p(\bar{t}+T) < p(\bar{t})$.
Similarly, for all $n \in \N^*$, we can show by induction $p(\bar{t} + nT) < p\bigl(\bar{t} + (n-1)T\bigr)$.
The sequence $\Bigl(p\bigl(\bar{t} + nT\bigr)\Bigr)_{n\in\N}$ is decreasing and bounded, thus it converges to a certain limit $\underline{p^\infty}$, which is a fixed point of $P$.

By uniqueness in Cauchy-Lipschitz theorem, the solution $p$ cannot cross $p^T_m$, hence, $\underline{p^\infty} = 0$.
The same arguments allow us to conclude that $\ds \lim_{t\to +\infty} p(t) = 0$.
%Otherwise, if $p_0 \in \left(p^T_m(\bar{t}), p^T_M(\bar{t})\right)$, then, according to the sign of $\Phi(p_0)$, one can use similar arguments as in \ref{1ercas} or \ref{2ecas}. If $\Phi(p_0) > 0$, then the sequence $\left(p(\bar{t}+nT)\right)_{n \in \N}$ is increasing and bounded from above, thus converges towards the closest zero of $\Phi$; if $\Phi(p_0) < 0$, then the sequence $\left(p(\bar{t}+nT)\right)_{n \in \N}$ is decreasing and bounded from below, thus converges towards the closest zero of $\Phi$. We deduce from these arguments that $p$ cannot converge, otherwise each sequence $\left(p(t_i+nT)\right)_{n \in \N}$, for $i$ varying, would converge towards the same limit, then there would be a constant solution to the ODE \eqref{ODE} in between $p_m^T$ and $p_M^T$. 
\end{enumerate}
\end{proof}

\subsection{Periodic solution for a nonlinearity with separated variables, proof of Proposition \ref{separation_var_unicite}} \label{periodic_solutions}

\begin{proof}[Proof of Proposition \ref{separation_var_unicite}]
Let $f : \R_+ \times [0,1] \to \R$ be defined by $f(t,p) := m(t) \tilde{f}(p)$, with $m$ and $\tilde{f}$ satisfying the hypotheses (H\ref{hyp_periodicite_bis}) to (H\ref{hyp_bistable_bis}).

First, note that the solution constantly equal to $\bar{\theta}$ is $T$-periodic and strictly between 0 and 1.
Let us work by contradiction and denote $p^T$ another $T$-periodic solution strictly between 0 and 1. By uniqueness in Cauchy-Lipschitz theorem, one has, for all $t \in \R_+, \; p^T(t) \notin \left\{ 0, \bar{\theta}, 1 \right\}$, which gives $\tilde{f}\left( p^T(t)\right) \neq 0$.%This implies that for all $t \in \R_+,\; (p^T)'(t) \neq 0$, which contradicts the fact that $p^T$ is $\C{1}$ and $T$-periodic.

One has for all $t \in \R_+, \; \ds \frac{(p^T)'(t)}{\tilde{f}\bigl(p^T(t)\bigr)} = m(t)$ therefore, integrating between 0 and $T$, one gets, denoting $F$ an antiderivative of $1 / \tilde{f}, \quad 0 = F\bigl(p^T(T)\bigr) - F\bigl(p^T(0)\bigr) = \int_0^T m(t) \dd t > 0$, which is a contradiction.

The non-attractivity of this constant solution can be inferred from the fact that $$\int_0^T\partial_p f(t,\bar{\theta}) \dd t = \tilde{f}'(\bar{\theta}) \int_0^T m(t) \dd t > 0$$ by (H\ref{hyp_periodicite_bis}) and (H\ref{hyp_bistable_bis}), using \cite{contri}, Proposition 2.3, which ends the proof of Proposition \ref{separation_var_unicite}.
\end{proof}

% Ne pas insérer de bibliographie à la fin de chaque sous-fichier

\section{Optimal control problem} \label{proofs_2}

\maketitle

The goal of this section is to find an asymptotic simplified version of the optimal control problem, that is equivalent to minimizing a seizable function, and to properly justify this convergence.

\subsection{Asymptotic control equation, proof of Proposition \ref{convergencep}} \label{proof_proposition_cv_p}

In this subsection, we study the asymptotic control equation when $\eps$ vanishes to 0.

First, let us formulate a quick remark on this proposition.
The positivity of $1 / g$ ensures that $p(t_0^+)$ is defined in a non-equivocal manner, given $p(t_0^-)$. Indeed, $G\bigl(p(t_0^+)\bigr) - G\bigl(p(t_0^{-})\bigr) = \ds \int_{p(t_0^-)}^{p(t_0^+)} \frac{1}{g}$, thus the choice of the antiderivative has no consequence.

\begin{proof}[Proof of Proposition \ref{convergencep}]
Let us first prove the pointwise convergence.

Set $t_f \in (t_0,+\infty)$. First, denote $(r_i)_{i \in \{0,...,n\}}$ a subdivision of $[0,t_0-\eps]$.

Let $i \in \{1,...,n\}$; for all $t \in [r_{i-1},r_i], \; \bigl|\ds (p^\eps)'(t)\bigr| = \bigl|f\bigl(t,p^\eps(t)\bigr)\bigr| \leqslant \norm{f}_\infty$.
Therefore $\bigl|p^\eps(r_i)-p^\eps(r_{i-1})\bigr| \leqslant \norm{f}_\infty |r_i - r_{i-1}|$, hence $\ds \sum_{i=1}^n \bigl|p^\eps(r_i)-p^\eps(r_{i-1})\bigr| \leqslant \norm{f}_\infty t_f $. This upper bound does not depend on the subdivision nor on $\eps$.

Denoting $(r'_i)_{i \in \{0,...,n\}}$ a subdivision of $[t_0,t_f]$, the same computations lead to the same conclusion. 

Finally, denote $(s_i)_{i \in \{0,...,n\}}$ a subdivision of the interval $[t_0-\eps,t_0]$.

Let $i \in \{1,...,n\}$. For all $t \in [s_{i-1},s_i], \; \bigl|\ds (p^\eps)'(t)\bigr| = \Bigl|f\bigl(t,p^\eps(t)\bigr) + \frac{C}{\eps K(t)} g\bigl(p^\eps(t)\bigr)\Bigr| \leqslant \norm{f}_\infty + \frac{C}{\eps \ds \min_{[0,T]} K}$.
Therefore, $\bigl|p^\eps(s_i)-p^\eps(s_{i-1})\bigr| \leqslant \biggl(\norm{f}_\infty + \ds \frac{C}{\eps \ds \min_{[0,T]} K}\biggr) |s_i - s_{i-1}|$.

Hence $\ds \sum_{i=1}^n \bigl|p^\eps(s_i)-p^\eps(s_{i-1})\bigr| \leqslant \biggl(\norm{f}_\infty + \ds \frac{C}{\eps \ds \min_{[0,T]} K}\biggr) \eps = \norm{f}_\infty \eps + \frac{C}{\ds \min_{[0,T]} K}$, which is bounded by an upper bound not depending on the subdivision nor on $\eps$.

Putting everything together, the sequence $(p^\eps)_\eps$ has a uniformly bounded total variation in $[0,t_f]$. Moreover, the sequence is uniformly bounded. Helly's selection theorem (see \cite{natanson}, section VIII.4; see \cite{belov-chistyakov}) therefore allows us to conclude that, in $[0,t_f]$, the sequence $(p^\eps)_\eps$ admits a pointwise convergent subsequence. This result being true for every interval $[0,t_f]$, with $t_f \in (t_0,+\infty)$, we can infer, by a diagonal argument, that in $[0,+\infty)$, the sequence $(p^\eps)_\eps$ admits a pointwise convergent subsequence, whose limit is denoted $p$.

Let us now investigate separately the two situations outside $t_0$ and at $t_0$ to get information on $p$.

\begin{description}

\item[Outside $\bm{t_0}$.]
For $\eps > 0$, by \eqref{controlsyst_g_wo_K_eps}, $p^\eps$ satisfies the equation $(p^\eps)' = f(t,p^\eps)$ in $(0,t_0-\eps) \cup(t_0,+\infty)$, alongside with the initial condition $p^\eps(0) = p_0$. %with $p'$ standing for the weak derivative of $p$.
Therefore, making $\eps$ go to 0, $p$ satisfies the initial condition $p(0) = p_0$ and the equation $p' = f(t,p)$ in $(0,t_0) \cup (t_0,+\infty)$ in the weak sense, with $p'$ standing for the weak derivative of $p$. Let us prove that this claim still holds in the classical sense.

Set $t_f \in (t_0,+\infty)$. Denoting $I := (0,t_0)$ (resp. $(t_0,t_f)$), one gets $\ds \int_I \ds |p'| = \ds \int_I \bigl| f(t,p) \bigr| < + \infty$ since $f$ is bounded. Therefore $p \in W^{1,1}(0,t_0)$ (resp. $W^{1,1}(t_0,t_f)$), $\ie$ $p$ is continuous in $(0,t_0) \cup (t_0,t_f)$ (rigorously it can be identified with the continuous representative of the equivalence class).
Hence, the function $t \mapsto f\bigl(t,p(t)\bigr)$ is continuous in $[0,t_0) \cup (t_0,t_f]$ and is equal to $p'$. This implies that $p'$ is continuous in $[0,t_0) \cup (t_0,t_f]$. This result being true for every $t_f \in (t_0,+\infty)$, we can infer that $p'$ is continuous in $[0,t_0) \cup (t_0,+\infty)$, therefore the equation $p' = f(t,p)$ is satisfied in the classical sense.

\item[At the time $\bm{t_0}$.]

%Let us first write the convergence of the sequence $(u^\eps)_\eps$.
%Setting $u := C \delta_{t_0}$, with $\delta_{t_0}$ the Dirac distribution at $t_0$,
%Let $\varphi \in \mathcal{D} := \C{\infty}_c([0,t_f])$. For all $\eps > 0, \; u^\eps \in L^1(0,t_f)$ hence it can be canonically identified with a distribution and:
%$$\langle u^\eps,\varphi \rangle _{\mathcal{D}',\mathcal{D}} - C \langle \delta_{t_0},\varphi \rangle _{\mathcal{D}',\mathcal{D}}  = \int_0^{t_f} u^{\eps}(t) \varphi(t) \dd t - C \varphi(t_0) = \frac{C}{\eps} \int_{t_0-\eps}^{t_0} \varphi(t) \dd t - C \varphi(t_0) \xrightarrow[\eps\to 0]{} 0.$$
%It is well known that the sequence $(u^\eps)_{\eps}$ weak-$*$ converges to $u$ on $\mathcal{D}'$ when $\eps \to 0$.

Set the change of variables $t = t_0 + (\tau - 1) \eps$.
Let us focus on the interval $t\in[t_0-\eps,t_0]$, $\ie$ $\tau \in [0,1]$. Set, for all $\tau \in [0,1], \; \tilde{p^\eps} (\tau) := p^\eps(t) = p^\eps \bigl(t_0 + (\tau-1)\eps\bigr)$.

One gets, for all $\tau \in [0,1]$,
$$(\tilde{p^\eps})'(\tau) = \eps f \bigl(t_0 + (\tau-1)\eps, \tilde{p^\eps}(\tau) \bigr) + \frac{C}{K\bigl(t_0+ (\tau-1)\eps\bigr)} \indic{[t_0-\eps,t_0]}\bigl(t_0+(\tau-1)\eps\bigr) g\bigl(\tilde{p^\eps}(\tau)\bigr).$$
The indicator function being constantly equal to 1, one gets, on $[0,1]$,
$$(\tilde{p^\eps})' = \eps f \bigl(t_0 + (\tau-1)\eps, \tilde{p^\eps} \bigr) + \frac{C}{K\bigl(t_0+ (\tau-1)\eps\bigr)}  g\bigl(\tilde{p^\eps}\bigr).$$

One can note that $\ds \bigl( (\tilde{p^\eps})' \bigr)_{\eps}$ is bounded since $f$ and $g$ are bounded and $K$ is bounded from below by a positive constant. Moreover, $\bigl(\tilde{p^\eps}\bigr)_{\eps}$ is bounded in $\C{0}\bigl([0,1]\bigr)$, thus $\bigl(\tilde{p^\eps}\bigr)_{\eps}$ is bounded in $\C{1}\bigl([0,1]\bigr)$. According to Arzelà-Ascoli theorem, it implies that $\bigl(\tilde{p^\eps}\bigr)_{\eps > 0}$ has a subsequence uniformly converging to some $\tilde{p}$. All those facts still hold for any subsequence of $\bigl(\tilde{p^\eps}\bigr)_{\eps > 0}$.

One also gets, for all $\tau \in [0,1]$,
$$\tilde{p^\eps}(\tau) = \tilde{p^\eps}(0) + \eps \int_0^\tau f \bigl(t_0 + (s-1)\eps, \tilde{p^\eps}(s) \bigr)\dd s + C \int_0^\tau \frac{g\bigl(\tilde{p^\eps}(s)\bigr)}{K\bigl(t_0+(s-1)\eps\bigr)}\dd s.$$

The boundedness of $f$ and the mean-value theorem applied to $g$ allow to pass to the limit, up to an extraction, when $\eps \to 0$ in this equation and to write, for all $\tau \in [0,1]$:

$$\tilde{p}(\tau)= \tilde{p}(0) + \frac{C}{K(t_0)} \int_0^\tau g\bigl(\tilde{p}(s)\bigr)\dd s,$$
$$\ie \quad \tilde{p}'(\tau) = \frac{C}{K(t_0)}g\bigl(\tilde{p}(\tau)\bigr).$$
Denoting $G$ an antiderivative of $1 / g$, one gets, in $[0,1]$, $(G \circ \tilde{p})' = \ds \frac{C}{K(t_0)}$.

Integrating between 0 and 1, it yields:
\begin{equation} \label{diff_Grondptilde}
G\bigl(\tilde{p}(1)\bigr) - G\bigl(\tilde{p}(0)\bigr) = \frac{C}{K(t_0)}.
\end{equation}
On the one hand, by pointwise convergence, one has, up to an extraction, $\tilde{p^\eps}(1) = p^\eps(t_0) \xrightarrow[\eps \to 0]{} p(t_0^+)$ and, up to an extraction, $\tilde{p^\eps}(1) \xrightarrow[\eps \to 0]{} \tilde{p}(1)$, therefore one concludes $\tilde{p}(1) = p(t_0^+)$ by double extraction.

On the other hand, let us prove that $\tilde{p}(0) = p(t_0^-)$.

One knows that the sequence $\bigl(\indic{[0,t_0-\eps]}\bigr)_{\eps}$ pointwise converges almost everywhere to $\indic{[0,t_0]}$, and by continuity of $f$ with respect to the second variable, one also gets that $f(t,p^\eps)$ pointwise converges to $f(t,p)$. Hence, by dominated convergence, it yields, up to an extraction:
$$\int_0^{t_0-\eps} f\bigl( t,p^\eps(t) \bigr) \dd t = 
\int_0^{t_0} \indic{[0,t_0-\eps]}(t) f\bigl( t,p^\eps(t) \bigr) \dd t \xrightarrow[\eps \to 0]{}
\int_0^{t_0} f\bigl( t,p(t) \bigr) \dd t.$$

Moreover, by \eqref{controlsyst_g_wo_K}, in $[0,t_0-\eps]$, one has $\ds (p^\eps)' = f(t,p^\eps)$, therefore, up to an extraction,
$$\tilde{p^\eps}(0) = p^\eps(t_0-\eps) = p^\eps(0) + \int_0^{t_0-\eps}f\bigl(t,p^\eps(t)\bigr) \dd t \xrightarrow[\eps \to 0]{} p(0) + \int_0^{t_0}f\bigl(t,p(t)\bigr) \dd t.$$
By uniqueness of the limit in the double extraction process, we obtain $\tilde{p}(0) = p(0) + \ds \int_0^{t_0}f\bigl(t,p(t)\bigr) \dd t$. Since in $[0,t_0)$, one has $p' = f(t,p)$, one gets $\tilde{p}(0) = p(t_0^-)$. Now one can rewrite \eqref{diff_Grondptilde} as follows, with a jump in the values of $G \circ p$:
$$ G\bigl(p(t_0^+)\bigr) - G\bigl(p(t_0^-)\bigr) = \frac{C}{K(t_0)}.$$
\end{description}
By both those conclusions, we infer that $p$ is unique, hence all pointwise convergent subsequences of $(p^\eps)_\eps$ have the same limit. For all $t \in (0,t_f)$, the sequence $\bigl(p^\eps (t) \bigr)_\eps$ has a unique accumulation point and the set $\{p^\eps(t), \eps >0\}$ has a compact closure; therefore the sequence $\bigl((p^\eps(t)\bigr)_\eps$ converges towards $p(t)$, that is $p^\eps$ pointwise converges towards $p$.
\end{proof}

\subsection{Convergence of the optimal control problem, proof of Theorem \ref{OCP_asymptotic}} \label{proof_theorem}

\begin{proof}[Proof of Theorem \ref{OCP_asymptotic}]

Thanks to Proposition \ref{basins_attraction}, one can rewrite $$\mathcal{I}_S = \Bigl\{C>0:\exists t_0^S \geqslant 1/S \;\text{ s. t. } p_S\text{ solution of \eqref{S_SCt0S} satisfies } 
p_S(t_0^S) > p^T_M(t_0^S) \Bigr\}.$$

\begin{description}
\item[Existence of the infimum.]
Here we prove that $\mathcal{I}_S$ is a non-empty set.
Let $\overline{K} = \ds \sup_{[0,T]} K$. We choose $\ds C > \frac{\overline{K}}{g\Bigl(\ds \max_{[0,T]}p^T_M\Bigr)}$. Let us prove that $C$ belongs to $\mathcal{I}_S$.

Set $t_0^S \geqslant \frac{1}{S}$. We want to prove that the solution of a simpler problem, which bounds from below the solution of our initial problem, is greater than $p^T_M$ at the time $t_0^S$.

One has $p^T_M\bigl(t_0^S\bigr) \leqslant \ds \max_{[0,T]} p^T_M < 1$. Indeed, $p^T_M$ does not reach 1, otherwise, by uniqueness in Cauchy-Lipschitz theorem, $p^T_M$ would be constantly equal to 1. Moreover, the function $g$ is decreasing, hence $\ds  0 < g\Bigl(\max_{[0,T]}p^T_M\Bigr) \leqslant \min_{[0,p^T_M(t_0^S)]} g$ for all $S$.
Therefore, one has $\ds C > \frac{\overline{K}}{\ds \min_{[0,p^T_M(t_0^S)]} g} \; p^T_M(t_0^S)$.

Let us denote $\underline{f}(p) := \ds \inf_{t\in[0,T]} f(t,p) < 0$. Let us take $S$ large enough such that $\ds C > \frac{\overline{K}}{\ds \min_{[0,p^T_M(t_0^S)]} g} \Biggl[ p^T_M(t_0^S) - \frac{\min_{[0,1]} \underline{f}}{S}\Biggr]$.

One has $\ds \frac{p^T_M(t_0^S)}{\min_{[0,1]} \underline{f} + \ds \frac{CS \ds \min_{[0,p^T_M(t_0^S)]} g}{\overline{K}}} < \frac{1}{S}$, which implies:
$$\frac{1}{S} > \ds \frac{p^T_M(t_0^S)}{\min_{[0,1]} \underline{f} + \ds \frac{CS \ds \min_{[0,p^T_M(t_0^S)]} g}{\overline{K}}} =
\bigints_0^{p^T_M(t_0^S)} \ds \frac{\dd s}{\min_{[0,1]} \underline{f} + \ds \frac{CS \ds \min_{[0,p^T_M(t_0^S)]} g}{\overline{K}}} \geqslant
\bigintss_0^{p^T_M(t_0^S)} \ds \frac{\dd s}{\ds \underline{f}(s) + \ds \frac{CS \ds g(s)}{\overline{K}}}.$$

Denoting now $\tilde{G}$ the antiderivative of $\ds \frac{1}{\underline{f} + \ds \frac{CSg}{\overline{K}}}$ vanishing at 0, it yields $\tilde{G} \Bigl(p^T_M\bigl(t_0^S\bigr)\Bigr) < \ds \frac{1}{S}$.

Let $\tilde{p}_S$ be the solution of the Cauchy problem
\begin{equation} \label{Stilde} \tag{$\mathcal{\tilde{S}}_{S,C,t_0^S}$}
	\begin{cases}
    \ds \tilde{p}_S' = \underline{f}(\tilde{p}_S) + \frac{C S}{\overline{K}} \indic{[t_0^S-\frac{1}{S},t_0^S]} \; g(\tilde{p}_S) \quad\text{a. e. in } (0,+\infty),\\
    \tilde{p}_S(0) = 0.
    \end{cases}
\end{equation}
One can note that $\tilde{p}_S$ exists and is unique since the function $(t,p) \mapsto \underline{f}(p) + \ds \frac{C S}{\overline{K}} \indic{[t_0^S-\frac{1}{S},t_0^S]}(t) \; g(p)$ is Lipschitz continuous with respect to the variable $p$ and locally integrable with respect to the variable $t$.
%Pour prouver que $\underline{f}$ est Lipschitz, revenir à la déf d'une fonction lipschitzienne et dire qu'il existe $t_2$ tq $\underline{f}(p_2) > f(t_2,p_2) - \eps$, puis majorer $\underline{f}(p_1)-\underline{f}(p_2)$

In the interval $[t_0^S-\frac{1}{S},t_0^S]$, \eqref{Stilde} gives $\ds \bigl(\tilde{G} \circ \tilde{p}_S\bigr)' = 1$.
Integrating this equality between $t_0^S-\frac{1}{S}$ and $t_0^S$ implies, since $\tilde{G}$ is increasing and vanishing at 0, $$\tilde{G} \Bigl(\tilde{p}_S\bigl(t_0^S\bigr)\Bigr) = \tilde{G} \Bigl(\tilde{p}_S\bigl(t_0^S\bigr)\Bigr) - \tilde{G} \biggl(\tilde{p}_S\Bigl(t_0^S-\frac{1}{S}\Bigr)\biggr) = \frac{1}{S} > \tilde{G}\Bigl(p^T_M\bigl(t_0^S\bigr)\Bigr).$$

Finally, as $\tilde{G}$ is invertible and increasing, it yields, for $S$ large enough, $\tilde{p}_S\bigl(t_0^S\bigr) > p^T_M\bigl(t_0^S\bigr)$.
Since $p_S$ is the solution of \eqref{S_SCt0S} and $\tilde{p}_S$ is the solution of \eqref{Stilde} and both problems have the same initial value, one has, for all $t>0,\; p_S(t) \geqslant \tilde{p}_S(t)$, hence 
$p_S\bigl(t_0^S\bigr) > p^T_M\bigl(t_0^S\bigr)$. $\mathcal{I}_S$ is a non-empty set therefore it has an infimum.

\item[Convergence of the sequence of infima to a certain $\bm{C_{\infty}}$.] Assuming that $S$ is large enough, we have proved that if $\ds C > \frac{\overline{K}}{g\Bigl(\ds \max_{[0,T]}p^T_M\Bigr)}$, then $C$ belongs to $\mathcal{I}_S$. This implies that $C_S^* \leqslant \ds \frac{\overline{K}}{g\Bigl(\ds \max_{[0,T]}p^T_M\Bigr)}$.

The sequence $(C_S^*)_S$ is bounded, therefore it converges, up to an extraction, when $S$ goes to $+\infty$. Let us denote $C_\infty$ this limit. We will prove that $C_{\infty} = C^*$ with a double inequality.

\item[Proof that $\bm{C_\infty \leqslant C^*}$.]
Let $\delta > 0$, and $p_S^\delta$ be the solution of the following control equation:
\begin{equation*}
	\begin{cases}
    \ds (p_S^\delta)'(t) = f\bigl(t,p_S^\delta(t)\bigr) + \frac{(C^* + \delta) S}{K(t)} \indic{[t_0^*-\frac{1}{S},t_0^*]}(t) \; g\bigl(p_S^\delta(t)\bigr) \quad \text{for almost every } t \in (0,+\infty),\\
    p_S^\delta(0) = 0.
    \end{cases}
\end{equation*}
By Proposition \ref{convergencep}, $\left(p_S^\delta\right)_S$ admits a pointwise convergent subsequence, whose limit, denoted $p^\delta$, is solution of the Cauchy problem
\begin{equation} \label{equationpinfdelta}
	\begin{cases}
    (p^\delta)'(t) = f\bigl(t,p^\delta(t)\bigr) \quad \text{for all } t \in (0,t_0^*) \cup (t_0^*,+\infty),\\
    G\bigl(p^\delta(t_0^*)\bigr)  = \ds \frac{C^*+\delta}{K(t_0^*)},\\
    p^\delta(0) = 0.
    \end{cases}
\end{equation}
This equation gives $G\bigl(p^\delta(t_0^*)\bigr) > \ds \frac{C^*}{K(t_0^*)} = G\bigl(p^T_M(t_0^*)\bigr)$ by definition of $C^*$, therefore $p^\delta(t_0^*) > p^T_M(t_0^*)$ since $G$ is increasing and invertible. Hence, there exists $S_\delta > 0$ such that for all $S \geqslant S_\delta, \; p_S^\delta(t_0^*) > p^T_M(t_0^*)$. This statement means that for all $S \geqslant S_\delta,\;C^* + \delta$ belongs to $\mathcal{I}_S$, therefore $C^* + \delta \geqslant C_S^*$.

Passing to the limit in the convergent subsequence yields $C^* + \delta \geqslant C_\infty$, this inequality being true for all $\delta > 0$; therefore $C^* \geqslant C_\infty$.

\item[Proof that $\bm{C_\infty \geqslant C^*}$.]
Assume that $C_\infty < C^*$. Then, there exists $\delta > 0$ such that $C_\infty < C^*(1-\delta)$, $\ie$ such that for a certain $S > 0$, one has $\inf \mathcal{I}_S = C_S^* < C^*(1-\delta)$. Hence there exists $C \in \mathcal{I}_S$ such that $C < C^*(1-\delta)$.

Let $t_0 \geqslant 1/S$ be such that the solution $p_S$ of \eqref{S_SCt0S} satisfies $p_S(t_0) > p^T_M(t_0)$. As $G$ is increasing, one has $G\bigl(p_S(t_0)\bigr) > G\bigl(p^T_M(t_0)\bigr)$, therefore there exists $\rho > 0$ such that
\begin{equation} \label{Gdepm}
G\bigl(p_S(t_0)\bigr) > \ds \frac{1}{1-\rho} G\bigl(p^T_M(t_0)\bigr).
\end{equation}

We want to bound from above the right-hand side of the equation satisfied by $p_S$. For that purpose, let us find a bound from below for $K$ and a bound from above for $f$.

By (H\ref{hyp_K}), there exists $\delta$ such that if $|t-t_0| \leqslant \delta$, then $|K(t)-K(t_0)| \leqslant \rho K(t_0)$. Let $S \geqslant 1 / \delta$; for all $t \in [t_0-\frac{1}{S},t_0]$, one has $K(t) \geqslant K(t_0)(1- \rho)$.

%To begin with, a first order Taylor expansion gives, for all $t \in [t_0-\frac{1}{S},t_0]$, writing $t = t_0 - h$, when $h$ vanishes to $0^+$, $K(t) = K(t_0) \ds \biggl(1-h\frac{K'}{K}(t_0)\biggr) + o(h)$.
%
%One has $\ds K(t_0)\biggl(1-h\frac{K'}{K}(t_0)\biggr)\geqslant K(t_0)\biggl(1-\ds \frac{1}{S} \sup_{[0,T]} \Bigl| \frac{K'}{K} \Bigr| \biggr)$ therefore, for $S$ large enough, 
%\begin{equation} \label{minorationK}
%K(t) \geqslant K(t_0)\biggl(1-\ds \frac{1}{S} \ds \sup_{[0,T]} \Bigl| \frac{K'}{K} \Bigr| \biggr).
%\end{equation}
%
%Let us consider $S$ large enough such that \eqref{minorationK} is satisfied and $\ds \frac{1}{S} \sup_{[0,T]} \Bigl| \frac{K'}{K} \Bigr| < \rho$. One gets, for all $t \in [t_0-\frac{1}{S},t_0],\; K(t) \geqslant K(t_0)(1-\rho)$.

Moreover, defining, for all $p \in [0,1], \; \bar{f}(p) := \ds \max_{t\in[0,T]} f(t,p)$, one has, in $[t_0-\frac{1}{S},t_0]$,
$$\ds p_S' = f(t,p_S) + \frac{C S}{K(t)} g(p_S) \leqslant \bar{f}(p_S) + \frac{C S}{K(t_0)(1-\rho)} g(p_S).$$

Here let us state an important remark. Assume that there exists $C_{\max} > 0$ such that the solution $p_S$ of the equation
\begin{equation} \label{S_SCmaxt0} \tag{$\mathcal{S}_{S,C_{\max},t_0}$}
	\begin{cases}
    p_S'(t) = f\bigl(t,p_S(t)\bigr) + \ds \frac{C_{\max} S}{K(t)} \indic{[t_0-\frac{1}{S},t_0]}(t) \; g\bigl(p_S(t)\bigr) \quad \text{for almost every } t\in(0,+\infty),\\
    p_S(0) = 0,
    \end{cases}
\end{equation}
reaches 1, $\ie$ satisfies $p_S(t_0) = 1$.

In the interval $[t_0-\frac{1}{S},t_0]$, the function $p_S$ is a solution of the Cauchy problem
\begin{equation} \label{syst_entre2}
	\begin{cases}
    p_S'(t) = f\bigl(t,p_S(t)\bigr) + \ds \frac{C_{\max} S}{K(t)} \; g\bigl(p_S(t)\bigr), \\
    p_S(t_0) = 1.
    \end{cases}
\end{equation}
The function constantly equal to 1 solves \eqref{syst_entre2}, therefore, by uniqueness in Cauchy-Lipschitz theorem, $p_S \equiv 1$ in $[t_0-\frac{1}{S},t_0]$. This gives $p_S(t_0-\frac{1}{S}) = 1$, which contradicts $p_S \equiv 0$ in $[0,t_0-\frac{1}{S})$.

Therefore, $p_S < 1$ in $[0,+\infty)$ and $\ds \min_{[0,p_{S,\max}]} g > 0$, denoting $p_{S,T,\max} := \ds \max_{[0,T]} p_S$.

Set $\chi := C^*(1-\delta) - C > 0$ and let us consider $S$ large enough such that $S \geqslant \ds \frac{\sup_{[0,T]}|K'|}{\rho K(t_0)}$ and $\ds \frac{\ds \max_{[0,1]} \bar{f}}{S} \leqslant \ds \frac{\chi}{K(t_0)(1-\rho)} \ds \min_{[0,p_{S,T,\max}]} g$.

One gets in $[t_0-\frac{1}{S},t_0]$, $$p_S' \leqslant \ds \frac{(C+\chi)S}{K(t_0)(1-\rho)} g(p_S) = \ds \frac{C^*(1-\delta)S}{K(t_0)(1-\rho)} g(p_S),$$
which implies $\ds \frac{K(t_0)(1-\rho)}{C^*(1-\delta) S}\;(G\circ p_S)' \leqslant 1$. Integrating between $t_0-\frac{1}{S}$ and $t_0$, it yields, thanks to \eqref{Gdepm}, and given that $p_S(t_0-\frac{1}{S}) = 0$:
$$\frac{K(t_0)}{C^*(1-\delta)S}G\bigl(p^T_M(t_0)\bigr) < \frac{K(t_0)(1-\rho)}{C^*(1-\delta)S}G\bigl(p_S(t_0)\bigr) \leqslant \frac{1}{S}.$$

Finally, one gets
$$K(t_0) G\bigl(p^T_M(t_0)\bigr) < C^*(1-\delta) < C^* = \inf_{[0,T]} \Bigl[K(\cdot) G \bigl( p^T_M(\cdot)\bigr)\Bigr],$$
which is a contradiction.

As the limit $C^*$ is uniquely determined, the entire sequence $(C^*_S)_S$ converges towards $C^*$, rather than only a subsequence.

\item[Convergence of the sequence of minimizers.]
$C_S^*$ is the infimum of $\mathcal{I}_S$, therefore there exists $t_0^{S*} \in [1/S,+\infty)$ such that the solution $p_S^*$ of $(\mathcal{S}_{S,C_S^*,t_0^{S*}})$ satisfies $\forall t<t_0^{S*}, \; p_S^*(t) < p^T_M(t)$, and $\forall t \geqslant t_0^{S*}, \;p_S^*(t) = p^T_M(t)$.

Indeed, let us work by contradiction and assume that for all $t_0^{S*} \in [1/S,+\infty)$, the solution $p_S^*$ of $(\mathcal{S}_{S,C_S^*,t_0^{S*}})$ either satisfies (a) $\exists t < t_0^{S*}:p_S^*(t) \geqslant p_M^T(t)$ or (b) $\exists t \geqslant t_0^{S*}:p_S^*(t) \neq p_M^T(t)$.
\begin{description}
\item[Assume (a).] As the control operates up to the time $t_0^{S*}$, one gets $p_S^*(t_0^{S*}) > p_M^T(t_0^{S*})$. Therefore $C_S^* \in \mathcal{I}_S$, and $\mathcal{I}_S$ is an open set, which is a contradiction.

\item[Assume (b).] It implies, by uniqueness in Cauchy-Lipschitz theorem, that (b1) for all $t \geqslant t_0^{S*}, \; p_S^*(t) < p_M^T(t)$ or (b2) for all $t \geqslant t_0^{S*}, \; p_S^*(t) > p_M^T(t)$.
(b1) implies the existence of a $C > C_S^*$ that bounds $\mathcal{I}_S$ from below, and (b2) gives $C_S^* \in \mathcal{I}_S$. Both are contradictions.
\end{description}

The sequence of minimizers $\bigl(t_0^{S*}\bigr)_S$ is bounded, therefore it converges up to an extraction to a limit, denoted $t_0^{\infty}$.
The goal is to bound $\ds \frac{C^*}{K(t_0^\infty)}$ by two quantities close enough to $G\bigl(p^T_M(t_0^\infty)\bigr)$.

For all $(t,p) \in [0,T] \times [0,1], \; \biggl| \ds \frac{f(t,p)K(t)}{p} \biggr| \leqslant \norm{\partial_p f}_{\infty} \ds \max_{[0,T]} K := M_f$.

Since $\ds (p_S^*)' = f(t,p_S^*) + \frac{C_S^* S}{K} \indic{[t_0^{S*}-\frac{1}{S},t_0^{S*}]} \; g(p_S^*)$ in $\R_+$, this bound gives, for all $t \geqslant 0$,
$$- \frac{M_f p_S^*(t)}{K(t)} + \frac{C_S^* S}{K(t)} \indic{[t_0^{S*}-\frac{1}{S},t_0^{S*}]}(t) \; g\bigl(p_S^*(t)\bigr)
\leqslant
(p_S^*)'(t)
\leqslant
\frac{M_f p_S^*(t)}{K(t)} + \frac{C_S^* S}{K(t)} \indic{[t_0^{S*}-\frac{1}{S},t_0^{S*}]}(t) \; g\bigl(p_S^*(t)\bigr).$$

One has $p_S^*(t_0^{S*}) = p_M^T(t_0^{S*}) < 1$, thus for all $t \in [t_0^{S*}-\frac{1}{S},t_0^{S*}]$, one gets $g\bigl(p_S^*(t)\bigr) \geqslant \ds \min_{[0,p_S^*(t_0^{S*})]} g = g\bigl(p_S^*(t_0^{S*})\bigr) > 0$.
In $[t_0^{S*}-\frac{1}{S},t_0^{S*}]$, as soon as $S$ is large enough to get
$C_S^* \geqslant \frac{C^*}{2}$ and $\ds -\frac{M_f}{K(t)} + \frac{C^*S}{2K(t)} \ds \min_{[0,p_S^*(t_0^{S*})]} g > 0$, one can write $\ds - \frac{M_f p_S^*(t)}{K(t)} + \frac{C_S^* S}{K(t)} \; g\bigl(p_S^*(t)\bigr) > 0$, therefore:
\begin{equation} \label{encad_composee}
\frac{(p_S^*)'(t)}{\ds \frac{M_f p_S^*(t)}{C_S^* S} + g\bigl(p_S^*(t)\bigr)}
\leqslant
\frac{C_S^* S}{K(t)}
\leqslant
\frac{(p_S^*)'(t)}{-\ds \frac{M_f p_S^*(t)}{C_S^* S} + g\bigl(p_S^*(t)\bigr)}.
\end{equation}
Setting, for all $p \in \bigl[0,\ds \max_{[0,T]} p^T_M\bigr], \; \underline{G_S}(p) := \ds \int_0^p \frac{\dd u}{\ds \frac{M_f u}{C_S^* S} + g(u)}$ and $\overline{G_S}(p) := \ds \int_0^p \frac{\dd u}{-\ds \frac{M_f u}{C_S^* S} + g(u)}$, the double inequality \eqref{encad_composee} reads:
$$\Bigl(\underline{G_S} \circ p_S^*\Bigr)'(t) 
\leqslant
\frac{C_S^* S}{K(t)}
\leqslant
\Bigl(\overline{G_S} \circ p_S^*\Bigr)'(t),$$

that we integrate between $t_0^{S*}-\frac{1}{S}$ and $t_0^{S*}$, getting, since $p_S^*$ is constantly equal to 0 until $t_0^{S*}-\frac{1}{S}$,
\begin{equation} \label{inegalite_GSbar}
\underline{G_S}\bigl(p_S^*(t_0^{S*})\bigr)
\leqslant C_S^* S \int_{t_0^{S*}-\frac{1}{S}}^{t_0^{S*}} \frac{\dd t}{K(t)}
\leqslant \overline{G_S}\bigl(p_S^*(t_0^{S*})\bigr).
\end{equation}
Now we have two bounds that seem close to $G\bigl(p^T_M(t_0^\infty)\bigr)$. The next step is to detail rigorously this approximation.

For all $u \in \bigl[0,\ds \max_{[0,T]} p^T_M\bigr]$, denote $\underline{g_S}(u) :=\ds \frac{M_f u}{C_S^* S} + g(u)$ and $\overline{g_S}(u) := \ds -\frac{M_f u}{C_S^* S} + g(u)$.
On the one hand, one has $\left(\underline{g_S}\right)_{S>0}$ uniformly converges to $g$ and $\ds \inf_{[0,\max_{[0,T]} p^T_M]} g > 0$, hence $\biggl(\ds \frac{1}{\underline{g_S}}\biggr)_{S>0}$ uniformly converges to $1/g$.
On the other hand, one has $\left(\overline{g_S}\right)_{S>0}$ uniformly converges to $g$ and $\overline{g_S} \geqslant - \ds \frac{M_f \ds \max_{[0,T]} p^T_M}{C_S^* S} + g$ and for $S$ large enough,
$\ds \inf_{[0,\max_{[0,T]}p^T_M]} \biggl[ - \ds \frac{M_f \ds \max_{[0,T]} p^T_M}{C_S^* S} + g \biggr] > 0$, hence $\biggl(\ds \frac{1}{\overline{g_S}}\biggr)_{S>0}$ uniformly converges to $\ds \frac{1}{g}$.

%Ces deux passages à l'inverse de convergence uniforme sont obtenus grâce à un lemme dans gy.tex

For all $p \in \bigl[0,\ds \max_{[0,T]} p^T_M\bigr]$, one has:
$$\Bigl|\underline{G_S}(p) - G(p)\Bigr| = \biggl|\int_0^p \biggl(\frac{1}{\underline{g_S}(u)} - \frac{1}{g(u)}\biggr) \dd u\biggr|
\leqslant
\int_0^p \norm{\frac{1}{\underline{g_S}} - \frac{1}{g}}_{\infty} \dd u
\leqslant
\norm{\frac{1}{\underline{g_S}} - \frac{1}{g}}_{\infty},
$$
which gives that $\underline{G_S}$ uniformly converges to $G$; it can be proved similarly that $\overline{G_S}$ uniformly converges to $G$.

Let $\eps >0$. There exists $S$ large enough such that $\norm{\underline{G_S}-G}_{\infty} \leqslant \eps$ and $\norm{\overline{G_S}-G}_{\infty} \leqslant \eps$.
For $S$ large enough, the inequality \eqref{inegalite_GSbar} yields:
$$ G\bigl(p^T_M(t_0^{S*})\bigr) - \eps
= G\bigl(p_S^*(t_0^{S*})\bigr) - \eps
\leqslant C_S^* S \int_{t_0^{S*}-\frac{1}{S}}^{t_0^{S*}} \frac{\dd t}{K(t)}
\leqslant G\bigl(p_S^*(t_0^{S*})\bigr) + \eps 
= G\bigl(p^T_M(t_0^{S*})\bigr) + \eps.$$

Thanks to the mean value theorem, there exists $\bar{t}_S \in (t_0^{S*}-\frac{1}{S},t_0^{S*})$ such that $S \ds \int_{t_0^{S*}-\frac{1}{S}}^{t_0^{S*}} \frac{\dd t}{K(t)} = \frac{1}{K(\bar{t}_S)}$.
Thus, by continuity of $K$, one has $C_S^* S \ds \int_{t_0^{S*}-\frac{1}{S}}^{t_0^{S*}} \frac{\dd t}{K(t)} \xrightarrow[S \to + \infty]{} \frac{C^*}{K(t_0^\infty)}$.

Moreover, thanks to the continuity of $p^T_M$ and $G$, one gets $G\bigl(p^T_M(t_0^{S*})\bigr) \xrightarrow[S \to + \infty]{} G\bigl(p^T_M(t_0^{\infty})\bigr)$, which finally yields: 
$$ G\bigl(p^T_M(t_0^{\infty})\bigr) - \eps
\leqslant \frac{C^*}{K(t_0^\infty)}
\leqslant G\bigl(p^T_M(t_0^{\infty})\bigr) + \eps.$$

Making $\eps$ vanish to 0, it gives $G\bigl(p^T_M(t_0^{\infty})\bigr) = \ds \frac{C^*}{K(t_0^\infty)}$, $\ie$ $t_0^\infty$ is a minimizer of the function $K(\cdot)G\bigl(p^T_M(\cdot)\bigr)$.
\end{description}
\end{proof}

\subsection{Nonlinearity with separated variables, proof of Proposition \ref{2releases}} \label{PFNL}

%Let us now study the particular case $f(t,p) = m(t) \tilde{f}(p)$, with $m$ and $\tilde{f}$ satisfying the hypotheses (H\ref{hyp_periodicite_bis}) to (H\ref{hyp_bistable_bis}).
%
%The unique $T$-periodic solution $p^T$ of \eqref{ODE} strictly between 0 and 1 being the function constantly equal to $\bar{\theta}$, by Proposition \ref{separation_var_unicite}, this allows us to rewrite the optimal control problem \eqref{P_eq} as follows:
%\begin{equation*}
%\ds \inf_{t_0 \in [0,T]} \Bigl[K(t_0) G(\bar{\theta}) \Bigr].
%\end{equation*}
%

\begin{proof}[Proof of Proposition \ref{2releases}]
Let us first prove the following lemma.
\begin{lemma} \label{decreasing_prop}
Let $0 < C_0 < \Cmin$ and $p$ be the solution of the control equation \eqref{Cauchy_IVP_0}, corresponding to a single release of $C_0$ individuals at $t_0$.

Then for all $t_1 > t_0,\; p(t_1) < p(t_0)$.
\end{lemma}

\begin{proof}

One has $C_0 < \ds \Cmin = G(\bar{\theta}) K(t_0^*)$, thus $\ds \frac{C_0}{K(t_0)} \leqslant \ds \frac{C_0}{K(t_0^*)} < G(\bar{\theta})$, which gives $\ds G^{-1}\Bigl(\frac{C_0}{K(t_0)}\Bigr) < \bar{\theta}$.

The equation \eqref{Cauchy_IVP_0} can be rewritten in $(t_0,+\infty)$, setting $p_0 := \ds G^{-1}\Bigl(\frac{C_0}{K(t_0)}\Bigr)$, as
\begin{equation*}
    \begin{cases}
    p'(t) = f\bigl(t,p(t)\bigr),\\
    p(t_0) = p_0 \text{ with } p_0 < \bar{\theta}.
    \end{cases}
\end{equation*}

Let us assume, by contradiction, that there exists a certain time $t_1 > t_0$ such that $p(t_1) = p(t_0)$.
Denoting $F$ an antiderivative of $\frac{1}{\tilde{f}}$, one has
$$0 = F\bigl(p(t_1)\bigr) - F\bigl(p(t_0)\bigr) = \ds \int_{t_0}^{t_1} p'(t) \frac{1}{\tilde{f}\bigl(p(t)\bigr)}\dd t = \ds \int_{t_0}^{t_1} m(t)\dd t > 0,$$

which is a contradiction.

Let $t_1 > t_0$. It has just been proved that $p(t_1) \neq p(t_0)$. But if $p(t_1) > p(t_0)$, then, since by Proposition \ref{basins_attraction}, one already knows that $\ds \lim_{t\to+\infty} p(t) = 0$, and since $p$ is continuous in $[t_0,+\infty)$, one gets the existence of a time $t_1' > t_1 > t_0$ such that $p(t_1') = p(t_0)$, which is a contradiction. Yet we conclude $p(t_1) < p(t_0)$.
\end{proof}

The first release $C_0$ at time $t_0$ generates a jump $G\bigl(p(t_0)\bigr)  = \ds \frac{C_0}{K(t_0)}$. The second release $C_1$ at time $t_1$ generates a jump
$G\bigl(p(t_1)\bigr) - G\bigl(p(t_1^-)\bigr)  = \ds \frac{C_1}{K(t_1)}$.
Summing those two jumps, one gets:
$$G\bigl(p(t_0)\bigr) - G\bigl(p(t_1^-)\bigr) + G\bigl(p(t_1)\bigr) =  \ds \frac{C_0}{K(t_0)} + \ds \frac{C_1}{K(t_1)}.$$
According to Lemma \ref{decreasing_prop}, one has $G\bigl(p(t_0)\bigr) > G\bigl(p(t_1^-)\bigr)$. Hence, $ G\bigl(p(t_1)\bigr)< \ds \frac{C_0}{K(t_0)} + \ds \frac{C_1}{K(t_1)} \leqslant \frac{\Cmin}{K(t_0^*)} = G(\bar{\theta})$, which yields $p(t_1) < \bar{\theta}$.
The system has not reached the threshold $\bar{\theta}$ at $t_1$, then, by Proposition \ref{basins_attraction}, it will converge to 0.
\end{proof}

% Ne pas insérer de bibliographie à la fin de chaque sous-fichier

\section{Application to population replacement with \textit{Wolbachia}} \label{application_wolbachia}

\maketitle

The mosquito \textit{Aedes aegypti} is one of the major vectors for transmitting dengue fever, chikungunya, Zika fever and yellow fever to the human being. None of those diseases can be fought with a specific treatment; moreover, there is no vaccine for Zika, and vaccines for chikungunya and dengue are only suitable for specific population groups (see \cite{OMSyellowfever}, \cite{OMSzika}, \cite{OMSdengue}, \cite{OMSchikungunya}). This is why an important part of prevention for those diseases consists in vector control.

Scientists have been studying the endosymbiotic bacterium \textit{Wolbachia}, mainly because of its interesting asset called \emph{pathogen interference} (PI), which consists in a significant decrease in the ability of \textit{Wolbachia}-infected mosquitoes to transmit the diseases when biting human beings (see \cite{bian}, \cite{moreira}, \cite{walker}). Therefore, replacing a wild mosquito population by a population of mosquitoes infected with \textit{Wolbachia} appears to be a valuable option to eradicate the transmission of those diseases; this \textit{Wolbachia}-based biocontrol is effective and does not seem detrimental on environment as insecticides can be. The World Mosquito Program has successfully implemented the \textit{Wolbachia} method in several countries in Asia and South America (see \cite{wmp}).

First, the bacterium is transmitted to mosquito eggs through micro-injections,
then it is maternally transmitted to the offspring. Furthermore, when an infected male mates with a non-infected female, the embryos die early in development, due to a phenomenon called \emph{cytoplasmic incompatibility} (CI) (see \cite{walker}). Consequently, if a population of infected mosquitoes is released among a wild population, it can be expected over time and under certain conditions, that, due to CI, almost all the mosquitoes will carry \textit{Wolbachia}: this is the \emph{population replacement} we are interested in.

It is therefore crucial to design release programs that pursue the goal of achieving the long-term persistence of \textit{Wolbachia}-infected mosquitoes, while optimizing economic resources associated with the breeding of infected mosquitoes and the releases.

Many mathematical works (see \cite{almeidachap}, \cite{almeida22}, \cite{almeida19}, \cite{orozco-gonzales}) have investigated optimization for \textit{Wolbachia}-based biocontrol, or population replacement, but they all considered the carrying capacity as a constant. However, the article \cite{dumont} deals with a time-varying carrying capacity; it underlines that the carrying capacity at a given time can be written as an affine function of the amount of available water at that time. In this section, we explain how we can, from the models given by population dynamics, get to a scalar control equation that meets with theoretical results proven in the previous sections.

Note that another well-known vector control strategy is the sterile insect technique and that some mathematical works (see \cite{bliman}) have explored optimal control approaches on this strategy.

\subsection{Scalar control equation}

For rigorous mathematical detail about this subsection, see \cite{almeida19} and \cite{strugarek}.
Computations are rewritten in a formal fashion, in order to retrieve the same results as in \cite{almeida19} and \cite{strugarek}, but here with a time-varying carrying capacity. Although those computations do not consist in a rigorous proof, they are sufficient to be convinced that with the same arguments as in those cited papers, we can obtain a similar model.
Moreover, the reader can note that the scalar model without control obtained in \cite{almeida19} exactly corresponds to the equation (6a) in \cite{barton}.

\label{scalar_control_equation}

\subsubsection{Competitive model}

To analyze mathematically the success of population replacement strategy, let us investigate here a control equation, where one acts on the wild population, releasing a time-distributed amount of infected individuals $u(t)$.

The evolution equations that we use model the competition of the released individuals with the wild ones. Let $n_1(t)$ denote the density of \textit{Wolbachia}-free mosquitoes (wild individuals) and $n_2(t)$ the density of \textit{Wolbachia}-infected mosquitoes at time $t$.
As introduced in \cite{almeida19}, mosquito population dynamics is modeled by the following competitive compartmental system:
\begin{equation} \label{mod_ini_1}
    \begin{cases}
    \ds n_1'(t) = b_1 n_1(t) \left(1-s_h\frac{n_2(t)}{n_1(t)+n_2(t)}\right) \left(1-\frac{n_1(t)+n_2(t)}{K(t)}\right) - d_1 n_1(t),\\[10pt]
    \ds n_2'(t) = b_2 n_2(t) \left(1-\frac{n_1(t)+n_2(t)}{K(t)}\right) - d_2 n_2(t) + u(t).
    \end{cases}
\end{equation}

In this system, $s_h$ denotes the CI rate, and one has $0 \leqslant s_h \leqslant 1$; $s_h = 1$ corresponding to a perfect CI, and $s_h = 0$ to no CI at all. For $i \in \{1,2\},\; b_i$ denotes the intrinsic birth rate, and $d_i$ the intrinsic mortality rate.

In this whole section, let us make an additional assumption on the function $K$:
\begin{hypothese} \label{hyp_Kreg}
$K$ is $\C{2}$ in $\R_+$.
\end{hypothese}

\subsubsection{Model reduction for large fecundity} \label{model-reduction}

\paragraph{Assumptions on the parameters.}
Let us detail several necessary assumptions on the parameters:

\begin{itemize}
\item Since \textit{Aedes} species has a very high reproductive power, it is relevant to consider that intrinsic birth rates are very large compared to intrinsic death rates, $\ie$ for $i \in \{1,2\}, \;b_i \gg d_i$, which makes us introduce $b_i^0$, comparable to $d_i$, be defined by $b_i = b_i^0/\sigma$, with $\sigma \ll 1$.

\item In the previous scaling, we always ensure that for $i\in\{1,2\},\; b_i^0 > d_i$, such that it fits the biological capacity of mosquito population to survive (see \cite{almeida19} and \cite{hughes}).

\item Mosquitoes infected by \textit{Wolbachia} have a slightly reduced fecundity and life span compared to wild mosquitoes (see \cite{hughes} and \cite{walker}), $\ie$ $b_1^0 > b_2^0$ and $d_1 < d_2$.

\item The CI rate $s_h$ satisfies the following condition (see \cite{almeida19}):
\begin{equation} \label{condition_sh}
1 - s_h < \frac{d_1b_2^0}{d_2b_1^0} < 1.
\end{equation}

This condition is biologically relevant. Indeed, \cite{dutra} and \cite{walker} show that CI is almost perfect since almost no egg hatch from a crossing between a non-infected female and an infected male.
\end{itemize}
%Note that we work on an interval $[0,t_f]$.

Let $n_1^\sigma$ and $n_2^\sigma$ be the solutions of \eqref{mod_ini_1}, with $b_1$ and $b_2$ given by the reduction depending on $\sigma$, and let $N^\sigma = n_1^\sigma + n_2^\sigma$ be the total amount of individuals and $\ds p^\sigma = \frac{n_2^\sigma}{N^\sigma}$ be the proportion of infected individuals, and $u$ the control. Denoting $\ds n^\sigma = \frac{1}{\sigma}\left(1-\frac{N^\sigma}{K}\right)$, system \eqref{mod_ini_1} reads: 
%\begin{equation*}
%    \begin{cases}
%    \ds
%    (n_1^\sigma)'(t) = \frac{b_1^0}{\sigma} n_1^\sigma(t) \bigl(1-s_h p^\sigma(t)\bigr) \left(1-\frac{N^\sigma(t)}{K(t)}\right) - d_1 n_1^\sigma(t),\\[10pt]
%    \ds (n_2^\sigma)'(t) = \frac{b_2^0}{\sigma} n_2^\sigma(t) \left(1-\frac{N^\sigma(t)}{K(t)}\right) - d_2 n_2^\sigma(t) + u(t).
%    \end{cases}
%\end{equation*}
\begin{equation} \label{mod_eps_2}
    \begin{cases}
    \ds
    (n_1^\sigma)'(t) = b_1^0 n_1^\sigma(t) \bigl(1-s_h p^\sigma(t)\bigr) n^\sigma(t) - d_1 n_1^\sigma(t),\\
    \ds
    (n_2^\sigma)'{t}(t) = b_2^0 n_2^\sigma(t) n^\sigma(t) - d_2 n_2^\sigma(t) + u(t).
    \end{cases}
\end{equation}

Summing both equations of \eqref{mod_eps_2}, replacing $n_2^\sigma$ by $p^\sigma N^\sigma$ on the one hand, and $n_1^\sigma$ by $\bigl(1-p^\sigma\bigr)N^\sigma$ on the other hand, and dividing by $N^\sigma$, one gets:
\begin{equation} \label{1surN_derivN}
    \frac{1}{N^\sigma}(N^\sigma)' = \Bigl[ b_1^0 \bigl( 1-p^\sigma \bigr) \bigl( 1-s_hp^\sigma \bigr) + b_2^0p^\sigma\Bigr]n^\sigma - \Bigl[ d_1\bigl( 1-p^\sigma \bigr) + d_2p^\sigma \Bigr] + \frac{u}{N^\sigma}.
\end{equation}
And since $(p^\sigma)' = \frac{1}{N^\sigma}(n_2^\sigma)' - p^\sigma \frac{1}{N^\sigma}(N^\sigma)'$, \eqref{mod_eps_2} implies:
\begin{equation} \label{approchee_derivp}
(p^\sigma)' = b_2^0 p^\sigma n^\sigma - d_2 p^\sigma + \frac{u}{N^\sigma} - p^\sigma  \frac{1}{N^\sigma}(N^\sigma)'.
\end{equation}

Moreover, one has $\ds n^\sigma = \frac{1}{\sigma}\left(1-\frac{N^\sigma}{K}\right)$, hence $\ds N^\sigma = K\bigl(1-\sigma n^\sigma\bigr)$ and $\ds \frac{1}{N^\sigma}(N^\sigma)' = \frac{1}{N^\sigma}K' \bigl( 1 - \sigma n^\sigma \bigr) -\sigma \frac{K}{N^\sigma} (n^\sigma)'$.

The following computations are formal, but are largely inspired by the work \cite{almeida19}, where all the convergences are proved rigorously when $K$ does not depend on time.
Assume that, when $\sigma$ vanishes to 0, $N^\sigma$ uniformly converges to $K$, $\ds \frac{1}{N^\sigma}(N^\sigma)'$ weakly converges to $\frac{K'}{K}$ up to an extraction.
Moreover, $p^\sigma$ uniformly converges to some $p$ when $\sigma$ vanishes to 0.

Furthermore, \eqref{1surN_derivN} gives:
\begin{equation*}
\ds
n^\sigma = \frac{
\ds\frac{1}{N^\sigma}(N^\sigma)' +  d_1\bigl( 1-p^\sigma \bigr) + d_2p^\sigma  - \frac{u}{N^\sigma}}
{b_1^0 \bigl( 1-p^\sigma \bigr) \bigl( 1-s_hp^\sigma \bigr) + b_2^0p^\sigma}.
\end{equation*}
The previous arguments give that $n^\sigma$ weakly converges to $\ds\frac{
\frac{K'}{K} +  d_1\bigl( 1-p \bigr) + d_2p  - \frac{u}{K}}
{b_1^0 \bigl( 1-p \bigr) \bigl( 1-s_hp \bigr) + b_2^0p}$.

Now, going to the limit in \eqref{approchee_derivp}, one gets:
\begin{equation*}
p' = b_2^0 p\frac{
\frac{K'}{K} +  d_1\bigl( 1-p \bigr) + d_2p  - \frac{u}{K}}
{b_1^0 \bigl( 1-p \bigr) \bigl( 1-s_hp \bigr) + b_2^0p} - d_2 p+ \frac{u}{K} - p  \frac{K'}{K},
\end{equation*}
\ie
\begin{equation*}
p' = p (1-p) \frac{
\frac{K'}{K}b_2^0 +  d_1 b_2^0 - d_2 b_1^0(1-s_h p) - \frac{1}{K} K' b_1^0 (1-s_h p)}
{b_1^0 \bigl( 1-p \bigr) \bigl( 1-s_h p \bigr) + b_2^0 p} 
+ \frac{u}{K} \frac{b_1^0 (1-p)(1-s_hp)}{b_1^0 \bigl( 1-p \bigr) \bigl( 1-s_h p \bigr) + b_2^0 p},
\end{equation*}
and $n_2(0) = 0$, thus $p(0) = 0$.

Let us define
\begin{equation*}
    \begin{cases}
    \ds
    f(t,p) = \frac{p(1-p)}{b_1^0(1-p)(1-s_h p) + b_2^0 p}\Biggl\{ s_h \Bigr(d_2b_1^0 + \frac{K'(t)}{K(t)}b_1^0\Bigr)p
- \biggl[ d_2b_1^0 - d_1b_2^0 + \frac{K'(t)}{K(t)}\bigl(b_1^0-b_2^0\bigr)\biggr]  \Biggr\},
    \\[10pt]
    \ds
    g(p) = \frac{b_1^0 (1-p)(1-s_hp)}{b_1^0 \bigl( 1-p \bigr) \bigl( 1-s_h p \bigr) + b_2^0 p}.
    \end{cases}
\end{equation*}
It gives:
\begin{equation*}
    \begin{cases}
    \ds
    p'(t) = f\bigl(t,p(t)\bigr) + \frac{u(t)}{K(t)} g\bigl(p(t)\bigr),\\
    p(0) = 0.
    \end{cases}
\end{equation*}

Denoting $\alpha(t) = s_h\Bigl(d_2 b_1^0 + \ds \frac{K'(t)}{K(t)} b_1^0\Bigr)$ and $\beta(t) = d_2 b_1^0 - d_1 b_2^0 + \ds \frac{K'(t)}{K(t)}(b_1^0-b_2^0)$, one gets:
\begin{equation} \label{deffetg}
	\begin{cases}
	\ds 
	f(t,p) = \frac{p(1-p)}{b_1^0(1-p)(1-s_h p) + b_2^0 p} \bigl( \alpha(t)p - \beta(t)\bigr),\\[10pt]
	\ds
	g(p) = \frac{b_1^0 (1-p)(1-s_hp)}{b_1^0 \bigl( 1-p \bigr) \bigl( 1-s_h p \bigr) + b_2^0 p}.
	\end{cases}
\end{equation}

Let us investigate whether those definitions ensure that $f$ satisfies the hypotheses (H\ref{hyp_regularite}) to (H\ref{hyp_vm_derivee_neg}) and (H\ref{hyp_bistable_fm}), and that $g$ satisfies (H\ref{hyp_g}).

By the expressions \eqref{deffetg}, it is straightforward that $f$ satisfies (H\ref{hyp_regularite}), (H\ref{hyp_periodicite}) and (H\ref{hyp_zeros}) and that $g$ satisfies (H\ref{hyp_g}).

Moreover, $\ds \partial_p f(t,0) = -\frac{\beta (t)}{b_1^0}$ gives:
$$\begin{array}{rcl}
\ds \frac{1}{T} \int_0^T \partial_p f(t,0) \dd t & = & - \ds \frac{1}{T} \int_0^T \frac{\beta(t)}{b_1^0} \dd t\\[10pt]
& = & - \ds \frac{1}{b_1^0 T} \int_0^T \biggl( d_2 b_1^0 - d_1 b_2^0 + \frac{K'(t)}{K(t)} \bigl(b_1^0 - b_2^0\bigr)\biggr) \dd t \quad \text{and $K$ is $T$-periodic}\\[10pt]
& = & \ds d_1 \frac{b_2^0}{b_1^0} - d_2 < 0 \quad \text{by \eqref{condition_sh}},
\end{array}$$

and $\ds \partial_p f(t,1) = \frac{\beta (t) - \alpha (t)}{b_2^0}$ gives:
$$\begin{array}{rcl}
\ds \frac{1}{T} \int_0^T \partial_p f(t,1) \dd t & = & \ds \frac{1}{T} \int_0^T \frac{\beta(t) - \alpha(t)}{b_2^0} \dd t\\[10pt]
& = & \ds \frac{1}{b_2^0 T} \int_0^T \biggl(d_2 b_1^0 - d_1 b_2^0 + \frac{K'(t)}{K(t)} \bigl(b_1^0 - b_2^0\bigr) -s_h \Bigl(d_2 b_1^0 + \frac{K'(t)}{K(t)}b_1^0\Bigr)\biggr) \dd t\\[10pt]
& = & \ds \frac{d_2b_1^0}{b_2^0} (1-s_h) - d_1 < 0 \quad \text{by \eqref{condition_sh}},
\end{array}$$

which implies that $f$ satisfies (H\ref{hyp_vm_derivee_neg}).

Knowing that $f_m(p) = \ds \frac{p(1-p)}{b_1^0(1-p)(1-s_hp)+b_2^0 p} \bigl(s_hd_2b_1^0p-d_2b_1^0+d_1b_2^0\bigr)$, we infer that $f_m$ vanishes at 0, 1 and $\bar{\theta} := \ds \frac{1}{s_h} \left( 1 - \ds \frac{d_1 b_2^0}{d_2 b_1^0} \right)$, and it is straightforward that $f_m < 0$ in $(0,\bar{\theta})$ and $f_m > 0$ in $(\bar{\theta},1)$. Moreover, (H\ref{hyp_vm_derivee_neg}) implies that $f_m'(0) < 0$ and $f_m'(1) < 0$, which yields by continuity $f_m'(\bar{\theta})>0$. The hypothesis (H\ref{hyp_bistable_fm}) is therefore satisfied.

Let us now mention the carrying capacity $K$. Considering that it strongly depends on the amount of available water in the environment (see \cite{dumont}), it seems relevant to consider that $K$ depends on the season of the year, this is why we assume from now that $K$ satisfies (H\ref{hyp_K}).

\subsection{Uniqueness of the periodic solution and minimization problem when the death rates are equal}

\begin{proposition} \label{death_rates_equal}
Let $f$ be the function defined by \eqref{deffetg} and assume that $d_1=d_2$. Then there exists a unique $T$-periodic solution strictly between 0 and 1 to the ODE \eqref{ODE}. Moreover, the optimal control problem \eqref{P_eq} is equivalent to minimizing the function $K$.
\end{proposition}

\begin{remark}
The case $d_1=d_2$ is not consistent with the biological assumptions, but is relevant for the study of the general case, that will be examined in the next subsection.
\end{remark}

\begin{proof}

The assumption $d_1 = d_2$ implies $\beta (t) = \bar{\theta} \alpha(t)$, with $\bar{\theta} = \ds \frac{b_1^0 - b_2^0}{s_hb_1^0}$, and the expression of $f$ becomes $$f(t,p) = \alpha(t) \ds \frac{p(1-p)}{b_1^0(1-p)(1-s_h p) + b_2^0 p} (p - \bar{\theta}) = \alpha (t) \tilde{f}(p),$$

with $\tilde{f}(p) := \ds \frac{p(1-p)}{b_1^0(1-p)(1-s_h p) + b_2^0 p} (p - \bar{\theta})$.

First, $K$ satisfies (H\ref{hyp_Kreg}) and $\int_0^T \alpha = s_h b_1^0 d_2 T > 0$, hence $\alpha$ satisfies (H\ref{hyp_periodicite_bis}).
Furthermore, one has $f_m(p) = \ds \frac{1}{T} \int_0^T f(t,p) \dd t = \ds \frac{\tilde{f}(p)}{T} \int_0^T \alpha (t) \dd t = s_h b_1^0 d_2 \tilde{f}(p)$, which gives $f_m'(p) = s_h b_1^0 d_2 \tilde{f}'(p)$.

Therefore, $\tilde{f}$ satisfies (H\ref{hyp_regularite_bis}), and we infer from (H\ref{hyp_vm_derivee_neg}) and (H\ref{hyp_bistable_fm}) that $\tilde{f}$ satisfies (H\ref{hyp_vm_derivee_neg_bis}) and (H\ref{hyp_bistable_bis}).
Proposition \ref{separation_var_unicite} proves the uniqueness of the $T$-periodic solution strictly between 0 and 1.
Moreover, we are in the simplified framework detailed in Subsection \ref{PFNL}, therefore the optimization problem \eqref{P_eq} is equivalent to minimizing the function $K$.
\end{proof}

\subsection{Dependence of the periodic solution(s) on the death rates difference}

Let us denote, for all $t\in\R_+, \; \theta(t) := \ds \frac{\beta(t)}{\alpha(t)}$. For all $t\in \R_+, \; f(t,\cdot)$ vanishes at $\theta(t)$.

In the following lemma, we generalize the result of the previous subsection with a positive and small death rates difference, which appears to be relevant considering that in the literature (see \cite{almeida19}, \cite{orozco-gonzales}), the relative difference $\frac{d_2-d_1}{d_2}$ is approximately equal to 10~\%.

\begin{lemma} \label{ecart_solper_thetabar}
Assume that the function $K$ satisfies (H\ref{hyp_Kreg}). Set $\eta = d_2 - d_1 > 0$.

For $\eta$ small enough fixed, the ODE \eqref{ODE} has a unique $T$-periodic solution strictly between 0 and 1, denoted $p^T_\eta$.

Let $\bar{\theta}_d = \ds \frac{b_1^0 - b_2^0}{s_hb_1^0}$ be the value of $\theta(t)$ in the limit case $d_1 = d_2$.
Then, one has for all $t \in [0,T], \; p_\eta^T(t) \underset{\eta \to 0}{=} \bar{\theta}_d + O(\eta)$.
\end{lemma}

\begin{proof}
For all $t \in \R_+$, one has
$$
\beta(t) = d_2 b_1^0 - (d_2-\eta)b_2^0 + \ds \frac{K'(t)}{K(t)}(b_1^0-b_2^0)
= \Bigl(d_2 + \ds \frac{K'(t)}{K(t)}\Bigr)(b_1^0-b_2^0) + \eta b_2^0
= \bar{\theta}_d \alpha(t) + \eta b_2^0.
$$

Replacing $\beta (t)$ by its expression, one gets, for all $(t,p) \in \R_+ \times [0,1]$,
$$
\begin{array}{rcl}
f(t,p)	& = & \ds \frac{p(1-p)}{b_1^0(1-p)(1-s_h p) + b_2^0 p} \bigl( \alpha(t)(p - \bar{\theta}_d) - \eta b_2^0 \bigr)\\[10pt]
		& = & \ds \alpha(t) \ds \frac{p(1-p)}{b_1^0(1-p)(1-s_h p) + b_2^0 p} (p - \bar{\theta}_d) - \eta b_2^0 \frac{p(1-p)}{b_1^0(1-p)(1-s_h p) + b_2^0 p}.
\end{array}
$$

Setting, for all $(t,p) \in \R_+ \times [0,1], \; f^0(t,p) := \alpha(t) \ds \frac{p(1-p)}{b_1^0(1-p)(1-s_h p) + b_2^0 p} (p - \bar{\theta}_d)$, $f^1(t,p) = f^1(p) := - b_2^0 \ds \frac{p(1-p)}{b_1^0(1-p)(1-s_h p) + b_2^0 p}$ and $f^\eta := f$, it gives $f^\eta = f^0 + \eta f^1$.

\begin{description}
\item[Uniqueness of the periodic solution.]
Denoting $\tilde{f}^0(p) := \ds \frac{p(1-p)}{b_1^0(1-p)(1-s_h p) + b_2^0 p} (p - \bar{\theta}_d)$, one has $f^0(t,p) = \alpha(t) \tilde{f}^0(p)$. One already knows that $p \equiv \bar{\theta}_d$ is a $T$-periodic solution of the ODE $p' = f^0(t,p)$.

Let us investigate its uniqueness. The function $\alpha(\cdot)$ satisfies (H\ref{hyp_periodicite_bis}). Moreover, straightforward computations yield $(\tilde{f}^0)'(0) < 0$, \;
$(\tilde{f}^0)'(1) < 0,\;$ and
$(\tilde{f}^0)'(\bar{\theta}_d) > 0$, therefore $\tilde{f}^0$ satisfies hypotheses (H\ref{hyp_regularite_bis}) to (H\ref{hyp_bistable_bis}); hence, by Proposition \ref{separation_var_unicite}, $p \equiv \bar{\theta}_d$ is the unique $T$-periodic solution strictly between 0 and 1 of the ODE $p' = f^0(t,p)$. Moreover, $\int_0^T \partial_p f^0\left(t,\bar{\theta}_d\right) \dd t > 0$ (see proof of Proposition \ref{separation_var_unicite}), that is the eigenvalue associated with this solution is negative.

Let us now verify the other conditions of Proposition \ref{perturbed}. 
Firstly, the function $f^0$, such as all the functions $f^\eta$, for $\eta > 0$, satisfy the hypotheses (H\ref{hyp_regularite}) to (H\ref{hyp_vm_derivee_neg}).
Furthermore, one has, for all $(t,p) \in \R_+ \times [0,1],\;\partial_p f^\eta(t,p) = \ds \partial_p f^0(t,p) + \eta \bigl(f^1\bigr)'(p)$.
Setting $\Delta := \norm{\bigl(f^1\bigr)'}_\infty$, which is finite since $\bigl(f^1\bigr)'$ is continuous on [0,1], it yields, for all $(t,p) \in \R_+ \times [0,1],\; \bigl| \partial_p f^0(t,p) - \partial_p f^\eta(t,p) \bigr| \leqslant \Delta \eta$. Denoting $\omega(\eta) := \Delta \eta$, one gets that $f^0$ and the family $\bigl(f^\eta\bigr)_{\eta>0}$ satisfy \eqref{condition_f_der_eps}.

Therefore, by Proposition \ref{perturbed}, for $\eta$ small enough, the ODE $p' = f^\eta(t,p)$ has a unique periodic solution between 0 and 1, which is not attractive.

\item[Speed of convergence.] Proposition 6.2 in \cite{contri} provides the existence of a constant $C_1$ such that for any $\eta > 0$, one has $\norm{P^0-P^\eta}_\infty < C_1 \Delta \eta$, where $P^0$ and $P^\eta$ respectively denote the Poincaré map associated with $f^0$ and $f^\eta$.
Setting $\Phi := P - \mathrm{Id}$, one immediately gets $\norm{\Phi^0-\Phi^\eta}_\infty < C_1 \Delta \eta$. 

Therefore, for all $\gamma >0$, one gets:
$$-C_1\Delta \eta < \Phi^0\left(\bar{\theta}_d +\gamma \eta\right) - \Phi^\eta\left(\bar{\theta}_d +\gamma\eta\right) < C_1 \Delta \eta,$$
hence
$$\Phi^\eta\left(\bar{\theta}_d +\gamma\eta\right) > \Phi^0\left(\bar{\theta}_d +\gamma\eta\right) - C_1 \Delta \eta.$$

The reader notes that the function $\Phi^0$ has the same regularity as the function $f^0$ by a generalization of Cauchy-Lipschitz theorem, therefore, in this application, $\Phi^0$ is at least $\C{2}\bigl([0,1]\bigr)$. Hence, for $\eta$ small enough, by Taylor remainder theorem, there exists $\theta_\eta^+ \in \left(\bar{\theta}_d,\bar{\theta}_d + \gamma \eta\right)$ such that: $$\Phi^\eta\left(\bar{\theta}_d +\gamma\eta\right) >
\Phi^0\left(\bar{\theta}_d\right) + \gamma \eta (\Phi^0)'\left(\bar{\theta}_d\right) + \frac{\gamma ^2 \eta ^2}{2} (\Phi^0)''(\theta_\eta^+)
- C_1 \Delta \eta 
=
\eta \left(\frac{\gamma ^2 \eta}{2} (\Phi^0)''(\theta_\eta^+) + \gamma (\Phi^0)'\left(\bar{\theta}_d\right)
- C_1 \Delta \right),
$$
since $p \equiv \bar{\theta}_d$ is an equilibrium state of $p' = f^0(t,p)$.

Set $\gamma > \ds \frac{C_1 \Delta}{(\Phi^0)'(\bar{\theta}_d)}$, this lower bound being positive since $(\Phi^0)'(\bar{\theta}_d) = \exp \left( \int_0^T \partial_p f^0(t,\bar{\theta}_d) \dd t \right) - 1 > 0$. 

If $(\Phi^0)''(\theta_\eta^+) \geqslant 0$, then for all $\eta > 0, \; \eta \left(\ds \frac{\gamma ^2 \eta}{2} (\Phi^0)''(\theta_\eta^+) + \gamma (\Phi^0)'\left(\bar{\theta}_d\right) - C_1 \Delta \right) > 0$. If $(\Phi^0)''(\theta_\eta^+) < 0$, let us take $0 < \eta < \ds \frac{\gamma (\Phi^0)'\left(\bar{\theta}_d\right) - C_1\Delta}{\frac{\gamma^2}{2} \max_{[\bar{\theta}_d,1]}\left|(\Phi^0)'' \right|}$,
it gives $0 < \eta < \ds \frac{C_1\Delta - \gamma (\Phi^0)'\left(\bar{\theta}_d\right)}{\frac{\gamma^2}{2} (\Phi^0)''(\theta_\eta^+)}$, and we infer $\eta \left(\ds \frac{\gamma ^2 \eta}{2} (\Phi^0)''(\theta_\eta^+) + \gamma (\Phi^0)'\left(\bar{\theta}_d\right) - C_1 \Delta \right) > 0$.

Finally, $\Phi^\eta\left(\bar{\theta}_d +\gamma\eta\right) > 0$ for $\eta$ small enough.

%It yields $$-C_1\Delta \eta < \Phi^0\left(\bar{\theta}_d +\ds \frac{C_1 \Delta}{(\Phi^0)'(\bar{\theta}_d)}\eta\right) - \Phi^\eta\left(\bar{\theta}_d +\ds \frac{C_1 \Delta}{(\Phi^0)'(\bar{\theta}_d)}\eta\right) < C_1 \Delta \eta,$$ therefore
%$$\Phi^\eta\left(\bar{\theta}_d +\ds \frac{C_1 \Delta}{(\Phi^0)'(\bar{\theta}_d)}\eta\right) > \Phi^0\left(\bar{\theta}_d +\ds \frac{C_1 \Delta}{(\Phi^0)'(\bar{\theta}_d)}\eta\right) - C_1 \Delta \eta.$$

%For $\eta$ small enough, a first order Taylor expansion gives: $$\Phi^\eta\left(\bar{\theta}_d +\ds \frac{C_1 \Delta}{(\Phi^0)'(\bar{\theta}_d)}\eta\right) > \Phi^0(\bar{\theta}_d) + \ds \frac{C_1 \Delta}{(\Phi^0)'(\bar{\theta}_d)} \eta (\Phi^0)'(\bar{\theta}_d)  - C_1 \Delta \eta = \Phi^0(\bar{\theta}_d) = 0$$ since $p \equiv \bar{\theta}_d$ is an equilibrium state of $p' = f^0(t,p)$.

Similarly, for $\eta$ small enough, by Taylor remainder theorem, there exists $\theta_\eta^- \in \left(\bar{\theta}_d - \gamma \eta,\bar{\theta}_d\right)$ such that: $$\Phi^\eta\left(\bar{\theta}_d -\gamma\eta\right) <
\eta \left(\frac{\gamma ^2 \eta}{2} (\Phi^0)''(\theta_\eta^-) - \gamma (\Phi^0)'\left(\bar{\theta}_d\right)
+ C_1 \Delta \right).
$$

If $(\Phi^0)''(\theta_\eta^-) \leqslant 0$, for all $\eta > 0, \; \eta \left(\ds \frac{\gamma ^2 \eta}{2} (\Phi^0)''(\theta_\eta^-) - \gamma (\Phi^0)'\left(\bar{\theta}_d\right) + C_1 \Delta \right) < 0.$
If $(\Phi^0)''(\theta_\eta^-) > 0$, taking $0 < \eta < \ds \frac{\gamma (\Phi^0)'\left(\bar{\theta}_d\right) - C_1\Delta}{\frac{\gamma^2}{2} \max_{[0,\bar{\theta}_d]}(\Phi^0)''}$, we get $\eta \left(\ds \frac{\gamma ^2 \eta}{2} (\Phi^0)''(\theta_\eta^-) - \gamma (\Phi^0)'\left(\bar{\theta}_d\right) + C_1 \Delta \right) < 0$.

Finally, $\Phi^\eta\left(\bar{\theta}_d -\gamma\eta\right) < 0$ for $\eta$ small enough.

We can now deduce, by continuity of $\Phi^\eta$, the existence of a real number $\theta_\eta \in \left(\bar{\theta}_d - \gamma \eta\;,\; \bar{\theta}_d + \gamma \eta  \right)$ such that $\Phi^\eta(\theta_\eta)=0$, $\ie$ such that the solution of $p'=f^\eta(t,p)$ starting at $\theta_\eta$ is a $T$-periodic solution, denoted $p_\eta^T$.

By the proof of Proposition 6.2 in \cite{contri}, one has, for all $t \in [0,T]$,
\begin{equation} \label{p0moinspT}
|p^0(t;\theta_\eta) - p^T_\eta(t)| = |p^0(t;\theta_\eta) - p^\eta(t;\theta_\eta)| \leqslant \ds C_1 \Delta \eta.
\end{equation}

Moreover, considering $\eta$ small enough and $t \in [0,T]$, one has
$$
\begin{array}{rcl}
\left|p^0(t; \theta_\eta) - p^0(t;\bar{\theta}_d)\right| 	& = 		& \left|\theta_\eta + \ds\int_0^t f^0\left(s,p^0(s;\theta_\eta)\right) \dd s - \bar{\theta}_d - \int_0^t f^0\left(s,p^0(s;\bar{\theta}_d)\right) \dd s\right|
\\
															& \leqslant & |\theta_\eta - \bar{\theta}_d| + \ds \int_0^t \left|f^0\left(s,p^0(s;\theta_\eta)\right) - f^0\left(s,p^0(s;\bar{\theta}_d)\right)\right| \dd s\\
															& \leqslant & |\theta_\eta - \bar{\theta}_d| + \ds \int_0^t \norm{\partial_p f^0}_\infty \left|p^0(t; \theta_\eta) - p^0(t;\bar{\theta}_d)\right| \dd s.
\end{array}
$$
by the mean value theorem.

Therefore, Grönwall lemma implies, for all $t \in [0,T]$,
\begin{equation} \label{p0diffDI}
\bigl|p^0(t; \theta_\eta) - p^0(t;\bar{\theta}_d)\bigr| \leqslant |\theta_\eta - \bar{\theta}_d| \exp \left(\norm{\partial_p f^0}_\infty t \right)  \leqslant \frac{C_1 \Delta}{(\Phi^0)'(\bar{\theta}_d)}\eta \exp \left(\norm{\partial_p f^0}_\infty T \right).
\end{equation}

Finally, one has, for $\eta$ small enough and for all $t \in [0,T]$,
$$ \bigl|p^T_\eta(t) - \bar{\theta}_d\bigr| \leqslant \bigl|p^T_\eta(t) - p^0(t; \theta_\eta)\bigr| + \bigl|p^0(t; \theta_\eta) - \bar{\theta}_d\bigr| = \bigl|p^0(t; \theta_\eta) - p^T_\eta(t)\bigr|  + \bigl|p^0(t; \theta_\eta) - p^0(t;\bar{\theta}_d)\bigr|,$$
that is, by \eqref{p0moinspT} and \eqref{p0diffDI},
$$\bigl|p^T_\eta(t) - \bar{\theta}_d\bigr| \leqslant C_1 \Delta \eta + \frac{C_1 \Delta}{(\Phi^0)'(\bar{\theta}_d)}\eta \exp \left(\norm{\partial_p f^0}_\infty T \right) = C_4 \Delta \eta,$$

where $C_4 := C_1\left(1 + \ds \frac{\exp \left(\norm{\partial_p f^0}_\infty T \right)}{(\Phi^0)'(\bar{\theta}_d)} \right)$.
\end{description}
\end{proof}

\subsection{Conclusions on the minimization problem in the general case}

In the general case, one has $d_1 = d_2-\eta$, with $\eta >0$. For a given $\eta$, let us denote $p^T_\eta$ the maximal $T$-periodic solution strictly between 0 and 1.

\begin{proposition} \label{etatendverszero}
Assume that the function $K$ satisfies (H\ref{hyp_Kreg}). For a given $\eta > 0$, denote $p^T_\eta$ the maximal $T$-periodic solution strictly between 0 and 1, and $t_\eta^*$ a minimizer of the function $K(\cdot)G\bigl(p_\eta^T(\cdot)\bigr)$. This minimizer satisfies the two following assertions:

\begin{enumerate}[label=(\roman*)]
\item \label{convergence_minimizer} The sequence $\left(t_\eta^*\right)_{\eta > 0}$ converges, up to an extraction, to a minimizer of $K$ in $[0,T]$, denoted $\tau_i^*$, when $\eta \to 0$.
\item \label{speed_convergence} Assume moreover that $K''$ is positive at each one of the minimizers of $K$. Up to the same extraction as in \ref{convergence_minimizer}, one has $|t_\eta^* - \tau_i^*| \underset{\eta \to 0}{=} O(\sqrt{\eta})$.
\end{enumerate}
\end{proposition}

The interpretation of \ref{convergence_minimizer} is the following: the best time to release the whole amount of infected
mosquitoes can become as close as we want to a time when the carrying capacity $K$ is minimal, as soon as $d_2 - d_1$ is small enough. A characterization of the speed of convergence is given in \ref{speed_convergence}.

\begin{proof} $\;$ %pour forcer le saut de ligne

\begin{enumerate}[label=(\roman*)]
\item The sequence of minimizers $\left(t^*_\eta\right)_{\eta>0}$ is bounded, therefore it converges, up to an extraction, towards a limit denoted $t_0^*$.
And Lemma \ref{ecart_solper_thetabar} implies that $p^T_\eta$ pointwise converges to $\thetabar$ when $\eta \to 0$.

For all $t \in [0,T]$, the inequality $K(t)G\bigl(p^T_\eta(t)\bigr) \geqslant K(t^*_\eta)G\bigl(p^T_\eta(t^*_\eta)\bigr)$ implies, making $\eta$ go to 0, \linebreak{$K(t)G(\thetabar) \geqslant K(t^*_0)G(\thetabar)$}, that is, since $G$ is positive, $t^*_0$ is a minimizer of $K$.

\item As $K''$ is positive at each one of its minimizers, each minimizer is isolated. Borel-Lebesgue theorem allows us to prove that $K$ has a finite number of minimizers.
%cf. gy.tex

Let us now denote $\left(\tau_i^*\right)_{i \in I}$ the finite family of minimizers of $K$.
Let us introduce the bounded function $q$ such that for all $t \in [0,T],\;p_\eta^T(t) = \thetabar + q(t) \eta$, and set $M_q := \norm{q}_\infty$.

Since for all $i \in I, \; K''(\tau_{i}^*) > 0$, by continuity, there exists $\zeta > 0$ such that for all $i \in I$, $K'' > 0$ in $[\tau_{i}^*-\zeta,\tau_{i}^*+\zeta]$. Denote $K''_\zeta := \ds \min_{i\in I} \min_{[\tau_{i}^*-\zeta,\tau_{i}^*+\zeta]} K'' > 0$

Set $C_5 > \ds 2 \ds\frac{M_q}{g(\thetabar)G(\thetabar)}\frac{\min_{[0,T]}K}{K''_\zeta}$. For $\eta > 0$, one has $$1+C_5 \eta \ds \frac{K''_\zeta}{\ds \min_{[0,T]}K} > 1 + 2\frac{M_q \eta}{g(\thetabar) G(\thetabar)}.$$

%Considering $t_i^*$ one of the accumulation points of the sequence $(t_\eta^*)_{\eta > 0}$,  denote $\mathcal{V}_{\eta,i}$ the following neighborhood of $t_i^*$:
%$$\mathcal{V}_{\eta,i} := \Biggl[t_i^*-2\sqrt{\frac{M_q }{g(\thetabar)G(\thetabar)} \frac{K(t_i^*)}{K''(t_i^*)} \eta}\quad,\quad t_i^*+2\sqrt{\frac{M_q }{g(\thetabar)G(\thetabar)} \frac{K(t_i^*)}{K''(t_i^*)} \eta}\Biggr],$$
%and let $t \notin \ds \bigcup_{i \in I} \mathcal{V}_{\eta,i}$, $\ie$ such that for all $i \in I, \;(t-t_i^*)^2 > \ds \frac{4 M_q}{g(\thetabar)G(\thetabar)} \frac{K(t_i^*)}{K''(t_i^*)} \eta$.

%It immediately infers that for all $i \in I, \; 1 +\ds \frac{(t-t_i^*)^2}{2} \frac{K''(t_i^*)}{K(t_i^*)} > 1 + 2\frac{M_q \eta}{g(\thetabar) G(\thetabar)}$.

Denoting, for $x>0, \; h(x) := \ds \frac{1+ \frac{M_q x}{g(\thetabar)}}{1- \frac{M_q x}{g(\thetabar)}}$, a first order Taylor expansion at 0 gives \linebreak{$h(x) \underset{x \to 0}{=} 1 + 2\ds \frac{M_q x}{g(\thetabar)}+o(x)$}, therefore for $\eta$ small enough, one then gets
$$1+C_5 \eta \ds \frac{K''_\zeta}{\ds \min_{[0,T]}K} > h\left(\frac{\eta}{G(\thetabar)}\right) = \ds \frac{G(\thetabar)+ \frac{M_q \eta}{g(\thetabar)}}{G(\thetabar) - \frac{M_q \eta}{g(\thetabar)}},$$ which implies $\ds \min_{[0,T]}K +C_5 \eta K''_\zeta > \min_{[0,T]}K \ds \frac{G(\thetabar)+ \frac{M_q \eta}{g(\thetabar)}}{G(\thetabar) - \frac{M_q \eta}{g(\thetabar)}}$.

%For $i \in I$ fixed, if $|t-t_i^*|$ is small enough, then, by a second order Taylor expansion, $K(t) \approx K(t_i^*) +\ds \frac{(t-t_i^*)^2}{2} K''(t_i^*)$ and $K(t) > K(t_i^*) \ds \frac{G(\thetabar)+ \frac{M_q \eta}{g(\thetabar)}}{G(\thetabar) - \frac{M_q \eta}{g(\thetabar)}}$, and this inequality still holds when $t$ is not close enough to $t_i^*$, even if it means choosing $\eta$ smaller, since $\ds \frac{G(\thetabar)+ \frac{M_q \eta}{g(\thetabar)}}{G(\thetabar) - \frac{M_q \eta}{g(\thetabar)}} \xrightarrow[\eta \to 0]{} 1$.

Hence, one has, for $\eta$ small enough, $$\Bigl(\min_{[0,T]}K +C_5 \eta K''_\zeta \Bigr) \left(G(\thetabar) - \frac{M_q \eta}{g(\thetabar)} \right) > \min_{[0,T]}K \left(G(\thetabar)+ \frac{M_q \eta}{g(\thetabar)}\right),$$
which gives, even if it means choosing $\eta$ smaller, by a first order Taylor expansion at 0,
\begin{equation} \label{1ere_etape}
\Bigl(\min_{[0,T]}K +C_5 \eta K''_\zeta \Bigr) \;G\bigl(\thetabar - M_q \eta\bigr) > \min_{[0,T]} K \times G\bigl(\thetabar + M_q \eta\bigr).
\end{equation}

Set $\nu := \sqrt{2C_5} > 0$.
Consider $\eta$ which makes \eqref{1ere_etape} true, and satisfying $\nu \sqrt{\eta} < \zeta$, and denote $\mathcal{V}_{\eta , i} := [\tau_i^*-\nu \sqrt{\eta},\tau_i^*+\nu \sqrt{\eta}]$.

Even if it means choosing $\eta$ smaller, for all $t \notin \ds \bigcup_{i \in I} \mathcal{V}_{\eta , i}$, one has $K(t) \geqslant \ds \min_{i\in I} \; \min \Bigl(K(\tau_i^*-\nu \sqrt{\eta}),K(\tau_i^*+\nu \sqrt{\eta})\Bigr)$, thus by Taylor remainder theorem, there exists $i_0 \in I$ and $\xi_{i_0} \in (\tau_{i_0}^*-\nu \sqrt{\eta},\tau_{i_0}^*+\nu \sqrt{\eta})$ such that
$$K(t) \geqslant K(\tau_{i_0}^*) + \ds \frac{\nu ^2 \eta}{2} K''(\xi_{i_0}) \geqslant K(\tau_{i_0}^*) + \ds \frac{\nu ^2 \eta}{2} K''_\zeta = \min_{[0,T]}K + C_5 \eta K''_\zeta.$$

Hence, by \eqref{1ere_etape}, one obtains, for all $t \notin \ds \bigcup_{i \in I} \mathcal{V}_{\eta , i}$:
$$K(t) G\bigl(\thetabar - M_q \eta\bigr) \geqslant \Bigl(\min_{[0,T]}K + C_5 \eta K''_\zeta \Bigr) \;G\bigl(\thetabar - M_q \eta\bigr) > \min_{[0,T]} K \times G\bigl(\thetabar + M_q \eta\bigr).$$

Furthermore, $\thetabar - M_q \eta \leqslant p_\eta^T(\cdot) \leqslant \thetabar + M_q \eta$ gives $G\bigl(\thetabar - M_q \eta\bigr) \leqslant G\bigl(p_\eta^T(\cdot)\bigr) \leqslant G\bigl(\thetabar + M_q \eta\bigr)$, which yields, for all $i \in I$,
$$K(t) G\bigl(p_\eta^T(t)\big) \geqslant K(t) \;G\bigl(\thetabar - M_q \eta\bigr) > \min_{[0,T]} K \times G\bigl(\thetabar + M_q \eta\bigr) \geqslant K(\tau_i^*) G\bigl(p_\eta^T(\tau_i^*)\bigr).$$

To  conclude, for $\eta$ small enough, one gets, for all $t \notin \ds \bigcup_{i \in I} \mathcal{V}_{\eta,i}$,
$K(t)G\bigl(p_\eta^T(t)\bigr) > K(\tau_i^*) G\bigl(p_\eta^T(\tau_i^*)\bigr)$.
Denoting $\tau^*$ the element of the family $\left(\tau_i^*\right)_{i \in I}$ such that $G\left(p_\eta^T(\tau^*)\right) = \ds \min_{i \in I} \Bigl[G\left(p_\eta^T(\tau_i^*)\right)\Bigr]$, one has, for $\eta$ small enough, and for all $t \notin \ds \bigcup_{i \in I} \mathcal{V}_{\eta,i}$, $K(t)G\bigl(p_\eta^T(t)\bigr) > K(\tau^*) G\bigl(p_\eta^T(\tau^*)\bigr)$.

Therefore, any minimizer $t_\eta^*$ of $K(\cdot)G\bigl(p_\eta^T(\cdot)\bigr)$ in $[0,T]$ belongs to $\ds \bigcup_{i \in I} \mathcal{V}_{\eta,i}$.
Finally, given a converging subsequence $\bigl(t_{\sigma(\eta)}^*\bigr)_{\eta>0}$, there exists $i \in I$ such that all the elements of the subsequence eventually belong to $\mathcal{V}_{\eta,i}$, which ends the proof.
%Given a converging subsequence $\bigl(t_{\sigma(\eta)}^*\bigr)_{\eta>0}$, there exists $\eta_0 \in \R_+^*$ such that there exists $i \in I$ such that for all $\eta < \eta_0$,  $t_{\sigma(\eta)}^* \in \mathcal{V}_{\eta,i}$, which ends the proof.
\end{enumerate}
\end{proof}

\subsection{Numerical simulations}
In this subsection, we implement an optimal control problem that aims to compute the control that minimizes its integral over time, remaining lower than a fixed upper bound, and allowing the solution to converge towards 1.
We use the package Gekko for Python (see \cite{gekko}), that relies on IPOPT (Interior Point OPTimizer) as a main solver for large-scale non-linear optimization problems, \textit{via} the optimization environment APMonitor.

\begin{description}

\item[Domain.]
Denoting $T$ the time period of the environmental carrying capacity, we set a time window of observation $[0,t_f]$ such that $t_f = NT$, with $N \geqslant 2$ an integer ($t_f$ could simply be greater than $T$).

\item[Numerical discretization.]
The control equation is implemented through an explicit Euler method,
%$\ie$, denoting $n$ the number of steps of length $h$ in the time interval $[0,t_f]$,
%$$\forall k \in \{0, ..., n-1\}, \quad t_k = kh, \quad p_{k+1} = p_k + h \left( f(t_k,p_k) + \frac{u_k}{K(t_k)} g(p_k) \right),\quad p_0 = 0,$$
%$p_k$ denoting an approximation of $p(t_k)$, with $p$ the solution of the control equation \eqref{controlsyst_g_wo_K}.
and the integral of $u$ is computed with a composite trapezoidal rule.
%$$\int_0^{t_f} u(t) \dd t \approx h \sum_{k=0}^{n-1}\frac{u_k + u_{k+1}}{2} = h \left(\sum_{k=0}^n u_k - \frac{1}{2}(u_0 + u_n)\right).$$

\item[Implemented functions.]
To make things simple, the carrying capacity, which is $T$-periodic and two times differentiable, is modeled as a cosine function:
$$K(t) = K_0 + A \cos \left(\frac{2\pi}{T}t\right),$$
with $K_0$ the mean value and $A$ the magnitude, such that $K_0 > A > 0$.

The implemented expressions of $f$ and $g$ are given by \eqref{deffetg}.

\item[Constraints.]
We set an upper bound for $u$, denoted $M$.
The control $u(\cdot)$ is set as constantly equal to 0 in $[T,t_f]$, to guarantee that the control occurs before $t=T$.
Moreover, in order to get a solution $p$ converging to 1, we implement the constraint $p(2T) > p(T)$. Let us prove that this condition is sufficient.

We consider a situation where \textbf{the $T$-periodic solution of the ODE \eqref{ODE} strictly between 0 and 1 is unique}, which is the case when the Poincaré map has the form plotted in Figure \ref{poincare}; let us denote it $p^T$.

As $u(\cdot) = 0$ in $[T,t_f]$, the solution $p$ of the control equation \eqref{controlsyst_g_wo_K} satisfies the ODE \eqref{ODE} in $[T,t_f]$.
Like in the proof of Proposition \ref{basins_attraction}, let us consider $P_T$ the Poincaré map shifted at time $T$, defined, for all $p_{0,T} \in [0,1]$, by $P_T(p_{0,T}) := p(T+T;p_{0,T})$, with $p(t ; p_{0,T})$ denoting, for all $t \in [T,t_f]$, the evaluation at $t$ of the solution $p$ of \eqref{ODE} in $[T,t_f]$ satisfying the condition $p(T) = p_{0,T}$.
Let us also denote $\Phi_T := P_T - \mathrm{Id}$.
Note that the condition $p(2T) > p(T)$ can be expressed as $\Phi_T\bigl(p(T)\bigr) > 0$.
By uniqueness of $p^T$, we infer that $p(T)$ belongs to the interval $(\theta_T,1)$, with $\theta_T := p^T(T)$, which gives, by Proposition \ref{basins_attraction}, that $p$ would converge to 1, if it was defined in $[0,+\infty)$.

\item[Objective.]
One wants to minimize the integral $\int_0^{t_f} u(t) \dd t$, given the aforementioned constraints.

\item[Results.]
The numerical value of the optimal control $u$ is plotted in Figure \ref{controleoptimal} for three different values of the upper bound: $M=0.02,\,0.04,\,0.06$. Moreover, the resulting value of $C$ is displayed, it is approximately equal to 0.01742, which appears, as expected, to be significantly close to 0.01739, the minimum of the function $K \times G \circ p^T$ numerically obtained by a gradient descent. In Figure \ref{solution}, the numerical approximation of the corresponding proportion $p$, solution of the control equation \eqref{controlsyst_g_wo_K} with the computed optimal control $u$, is plotted.
  
The optimal control $u$ obtained by this algorithm happens to be a short impulse with the form of a scaled indicator function showing a plateau at its upper bound $M$. It means that the optimal control is also the one of shortest duration. The instant of this optimal appears to be quite close to the minimizer of the function $K$ in the interval $[0,T]$, as established in Proposition \ref{etatendverszero}.
\begin{figure}[H] 
    \centering
    \includegraphics[height=8cm]{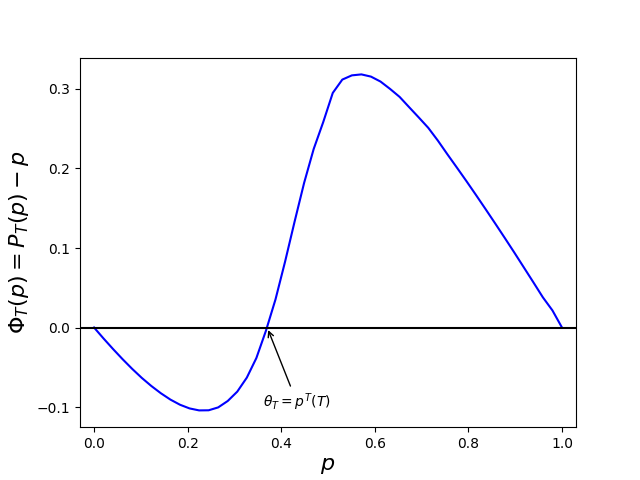}
    \caption{Graph of $\Phi_T = P_T - \mathrm{Id}$, with $P_T$ the Poincaré map shifted at time $T$, obtained by a Python script.
    %with the function \texttt{solve}\_\texttt{ivp} from the sub-package \texttt{integrate} of the SciPy library.
    The numerical values implemented in this script are $K_0 = 0.06, \; T = 12, \; b_1^0 = 0.8, \; b_2^0 = 0.6, \; d_1 = 0.27, \; d_2 = 0.3$ and $s_h = 0.9$, meeting with data found in the literature (see \cite{almeida19}, \cite{duprez}), and $A = 0.02$ arbitrarily. For this implementation, we take 50 points in the interval $[0,1]$.}
    \label{poincare}
\end{figure}

\begin{figure}[H] 
    \centering
    \includegraphics[height=8cm]{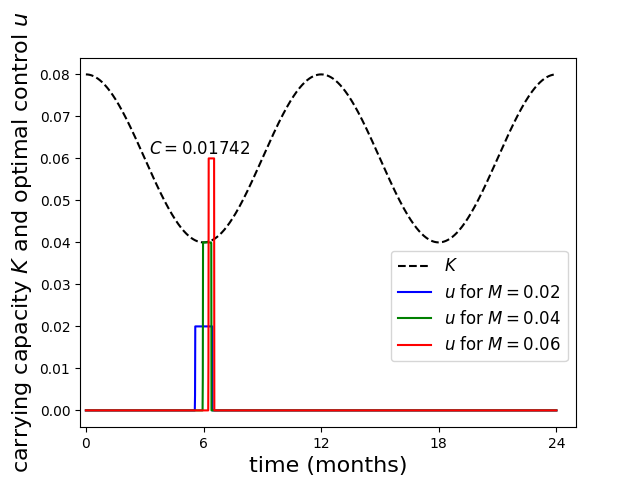}
    \caption{Optimal control $u$ obtained by a Python script with the package Gekko. The numerical values are exactly the same as in Figure \ref{poincare}. Moreover, we set $N = 2$ and $n = 1000$, and consider three different scenarii: $M = 0.02, M = 0.04$, and $M = 0.06$.}
    \label{controleoptimal}
\end{figure}

\begin{figure}[H]
    \centering
    \includegraphics[height=8cm]{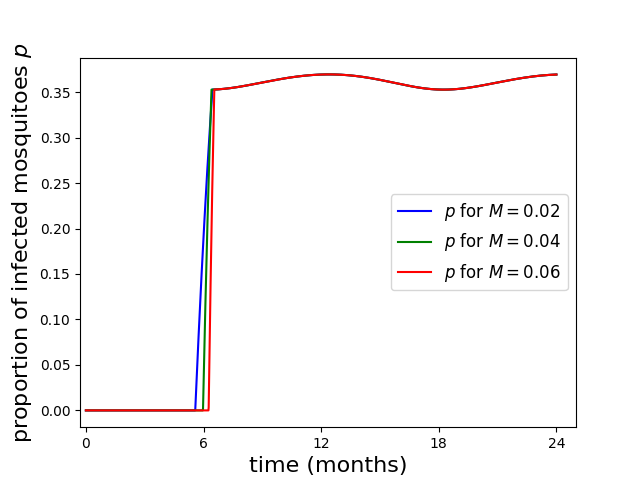}
    \caption{Optimal solution $p$, proportion of \textit{Wolbachia}-infected mosquitoes through time, obtained by the same Python script with Gekko. The numerical values are exactly the same as in Figure \ref{controleoptimal}.}
    \label{solution}
\end{figure}

In Subsection \ref{main_results-OCP}, in the theoretical search of optimal controls, we limited ourselves to control terms having the form of scaled indicator functions, and asymptotically Dirac impulses.
As far as population replacement with \textit{Wolbachia} is concerned, the relevance of this theoretical assumption is confirmed by numerical simulations.

\end{description}
% Ne pas insérer de bibliographie à la fin de chaque sous-fichier

\textbf{Acknowledgements.}

\begin{wrapfigure}{l}{0.15\textwidth}
\includegraphics[width=1\linewidth]{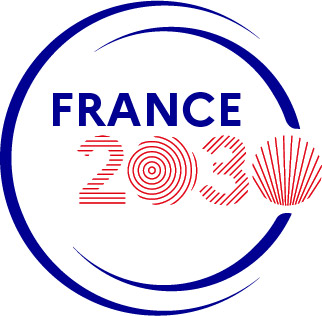} 
\end{wrapfigure}

G. N. was partially supported by the ANR project ReaCh (\textit{Réaction-diffusion : nouveaux Challenges}) and the ANR project Optiform (\textit{Optimisation de formes}), D. N. was supported by a grant from ANRS MIE (AAP 2024 - Fellowship Programme 2024 - FRAME) and N. V. was partially supported by the STIC AmSud project BIO-CIVIP (Biological control of insect vectors and insect pests) 23-STIC-02 and the ANR project Maths-ArboV (\textit{Modélisation mathématique et outils de gestion des
arboviroses émergentes}) 24-EXMA-0004, through the PEPR Maths-Vives of the program France 2030.

\printbibliography

@article{almeida19,
    author  = {L. Almeida and Y. Privat and M. Strugarek and N. Vauchelet},
    title   = {Optimal releases for population replacement strategies: \uppercase{a}pplication to \textit{Wolbachia}},
    year    = "2019",
    journal = "SIAM Journal on Mathematical Analysis",
    volume = "51-4",
    pages = "3170--3194"
}

@article{almeida22,
  author = {L. Almeida and J. Bellver Arnau and Y. Privat},
  title = "Optimal control strategies for bistable \uppercase{ode} equations: Application to mosquito population replacement",
  journal = "Applied Mathematics and Optimization",
  year = "2023",
  volume = "87",
  number = "10"
}

@article{almeida24,
    author  = {L. Almeida and J. Bellver Arnau and Y. Privat and C. Rebelo},
    title   = {Vector-borne disease outbreak control via instant releases},
    year    = "2024",
    journal = "Journal of Mathematical Biology",
    volume = "89-6",
    number = "63"
}

@inbook{almeidachap,
	title = {Minimal cost-time strategies for mosquito population replacement}, 
	booktitle = {Optimization and Control for Partial Differential Equations},
	author = {L. Almeida and J. Bellver Arnau and M. Duprez and Y. Privat},
	publisher = {De Gruyter},
	address = {Berlin, Boston},
	pages = {73--90},
	year = {2022},
	isbn={},
}

@article{contri,
    author  = {B. Contri},
    title   = {Pulsating fronts for bistable on average reaction-diffusion equations in a time periodic environment},
    year    = "2016",
    journal = "Journal of Mathematical Analysis and Applications",
    volume = "437-1",
    pages = "90--132"
}

@article{duprez,
	author = {M. Duprez and R. Hélie and Y. Privat and N. Vauchelet},
	title = {Optimization of spatial control strategies for population replacement, application to \textit{Wolbachia}},
	journal = {ESAIM: COCV},
	year = 2021,
	volume = 27,
	number = 74
}

@article{walker,
	author = {T. Walker and P. H. Johnson and L. A. Moreira and I. Iturbe-Ormaetxe and F. D. Frentiu},
	title = "The \textit{w}Mel \textit{Wolbachia} strain blocks dengue and
invades caged \textit{Aedes aegypti} populations",
	journal = {Nature},
	year = 2011,
	volume = 476,
	pages = "450--453",
}

@article{dutra,
    author = {H. L. C. Dutra and L. M. B. dos Santos and E. P. Caragata and J. B. L. Silva and D. A. M. Villela and R. Maciel-de-Freitas and L. A. Moreira},
    journal = {PLOS Neglected Tropical Diseases},
    title = {From lab to field: the influence of urban landscapes on the invasive potential of \textit{Wolbachia} in Brazilian \textit{Aedes aegypti} mosquitoes},
    year = {2015},
    volume = {9},
    pages = {1--22},
}

@article{hughes,
    author  = {H. Hughes and N. F. Britton},
    title   = {Modelling the use of \textit{Wolbachia} to control dengue fever transmission},
    year    = "2013",
    journal = "Bulletin of Mathematical Biology",
    volume = "75",
    pages = "796--818"
}

@book{natanson,
    title={Theory of functions of a real variable},
    author={I. P. Natanson},
    isbn={},
    series={New York},
    year={1965},
    publisher={Frederick Ungar Publishing Co.}
}

@article{belov-chistyakov,
    author  = {S. A. Belov and V. V. Chistyakov},
    title   = {A selection principle for mappings of bounded variation},
    year    = "2000",
    journal = "Journal of Mathematical Analysis and Applications",
    volume = "249-2",
    pages = "351--366"
}

@article{nedeljkov,
    author  = {M. Nedeljkov and M. Oberguggenberger},
    title   = {Ordinary differential equations with delta function terms},
    year    = "2012",
    journal = "Publications de l'Institut Mathématique. Nouvelle série",
    volume = "91--105",
    pages = "125--135"
}

@article{strugarek,
	author = {M. Strugarek and N. Vauchelet},
	title = {Reduction to a single closed equation for 2-by-2 reaction-diffusion systems of \uppercase{L}otka--\uppercase{V}olterra type},
	journal = "SIAM Journal on Applied Mathematics",
	volume = "76",
	number = "5",
	pages = "2060--2080",
	year = "2016",
}

@article{barton,
	author = {N. H. Barton and M. Turelli},
	title = {Spatial waves of advance with bistable dynamics: \uppercase{c}ytoplasmic and 	genetic analogues of \uppercase{A}llee effects},
	journal = "The American Naturalist",
	volume = "178",
	number = "3",
	pages = "E48--E75",
	year = "2011",
}

@article{orozco-gonzales,
	author = {J. L. Orozco Gonzales and A. dos Santos Benedito and H. de Oliveira Florentino and C. Pio Ferreira and D. Cardona-Salgado and L. S. Sepulveda-Salcedo and O. Vasilieva},
	title = {Optimization approaches to \textit{Wolbachia}-based biocontrol},
	journal = "Applied Mathematical Modelling",
	volume = "137",
	year = "2025",
}

@article{bian,
    author = {G. Bian and Y. Xu and P. Lu and Y. Xie and Z. Xi},
    journal = {PLOS Pathogens},
    publisher = {Public Library of Science},
    title = {The endosymbiotic bacterium \textit{Wolbachia} induces resistance to dengue virus in \textit{Aedes aegypti}},
    year = {2010},
    volume = {6},
    pages = {1-10},
    number = {4},
}

@article{moreira,
    author = {L. A. Moreira and I. Iturbe-Ormaetxe and J. A. Jeffery and G. Lu and A. T. Pyke and L. M. Hedges and B. C. Rocha and S. Hall-Mendelin and A. Day and M. Riegler and L. E. Hugo and K. N. Johnson and B. H. Kay and E. A. McGraw and A. F. van den Hurk and P. A. Ryan and S. L. O'Neill},
    journal = {Cell},
    publisher = {Public Library of Science},
    title = {A \textit{Wolbachia} symbiont in \textit{Aedes aegypti} limits infection with dengue, chikungunya, and plasmodium},
    year = {2009},
    volume = {139},
    pages = {1268--1278},
    number = {7},
}

@article{aronson-weinberger,
    author = {D. G. Aronson and H. F. Weinberger},
    journal = {Lecture Notes in Mathematics},
    publisher = {Springer},
    title = {Nonlinear diffusion in population genetics, combustion, and nerve pulse propagation},
    year = {1975},
    volume = {446},
    pages = {5--49},
    number = {1},
}

@article{alikakos-bates-chen,
    author = {N. D. Alikakos and P. W. Bates and X. Chen},
    journal = {Transactions of the American Mathematical Society},
    publisher = {},
    title = {Periodic traveling waves and locating oscillating patterns in multidimensional domains},
    year = {1999},
    volume = {351},
    pages = {2777--2805},
    number = {7},
}

@article{wang,
	title = {Speeds of invasion in a model with strong or weak Allee effects},
	journal = {Mathematical Biosciences},
	volume = {171},
	number = {1},
	pages = {83-97},
	year = {2001},
	author = {M.-H. Wang and M. Kot},
}

@article{gekko,
	title = {GEKKO Optimization Suite},
	author = {L. Beal and D. Hill and R. Martin and J. Hedengren},
	journal = {Processes},
	volume = {6},
	number = {8},
	pages = {106},
	year = {2018},
	publisher = {Multidisciplinary Digital Publishing Institute}
}

@article{schraiber,
	title = {Constraints on the use of lifespan-shortening \textit{Wolbachia} to control 	dengue fever},
	journal = {Journal of Theoretical Biology},
	volume = {297},
	pages = {26-32},
	year = {2012},
	author = {J. G. Schraiber and A. N. Kaczmarczyk and R. Kwok and M. Park and R. Silverstein and F. U. Rutaganira and T. Aggarwal and M. A. Schwemmer and C. L. Hom and R. K. Grosberg and S. J. Schreiber},
}

@misc{OMSdengue,
  author = {World Health Organization},
  title = {Dengue and severe dengue - Questions and answers},
  year = 2024,
  howpublished = {\url{https://www.who.int/news-room/questions-and-answers/item/dengue-and-severe-dengue}},
  note = {Accessed 2025/07/29}
}

@misc{OMSchikungunya,
  author = {World Health Organization},
  title = {Chikungunya},
  year = 2025,
  howpublished = {\url{https://www.who.int/news-room/fact-sheets/detail/chikungunya}},
  note = {Accessed 2025/07/29}
}

@misc{OMSzika,
  author = {World Health Organization},
  title = {Zika virus},
  year = 2022,
  howpublished = {\url{https://www.who.int/news-room/fact-sheets/detail/zika-virus}},
  note = {Accessed 2025/07/29}
}

@misc{OMSyellowfever,
  author = {World Health Organization},
  title = {Yellow fever - Questions and answers},
  year = 2017,
  howpublished = {\url{https://www.who.int/news-room/questions-and-answers/item/yellow-fever}},
  note = {Accessed 2025/07/29}
}

@article{dumont,
	author = {Y. Dumont and M. Duprez},
	title = {Modeling the impact of rainfall and temperature on sterile insect	control strategies in a tropical environment},
	journal = {Journal of Biological Systems},
	volume = {32},
	number = {01},
	pages = {311-347},
	year = {2024},
}

@article{almeida20,
	author = {L. Almeida and A. Haddon and C. Kermorvant and A. Léculier and Y. Privat and M. Strugarek and N. Vauchelet and J. P. Zubelli},
	title = {Optimal release of mosquitoes to control dengue transmission},
	journal = {ESAIM: Proceedings and Surveys},
	volume = {67},
	number = {},
	pages = {16-29},
	year = {2020},
}

@article{nadin,
  author = {G. Nadin and A. I. Toledo Marrero},
  title = "On the maximization problem for solutions of reaction-diffusion equations with respect to their initial data",
  journal = "Mathematical Modelling of Natural Phenomena",
  year = "2020",
  volume = "15",
  number = "71"
}

@article{mazari,
  author = {I. Mazari and G. Nadin and A. I. Toledo Marrero},
  title = "Optimisation of the total population size with respect to the initial condition for semilinear parabolic equations: two-scale expansions and symmetrisations",
  journal = "Nonlinearity",
  year = "2021",
  volume = "34",
  number = "11",
  pages = {7510-7539}
}

@misc{wmp,
  author = {World Moquito Program},
  title = {Our \textit{Wolbachia} method},
  howpublished = {\url{https://www.worldmosquitoprogram.org/en/work/wolbachia-method}},
  note = {Accessed 2025/11/18}
}

@article{pouchol,
	author = {C. Pouchol and E. Trélat and E. Zuazua},
	title = {Phase portrait control for 1D monostable and bistable reaction-diffusion equations},
	journal = {Nonlinearity},
	volume = {32},
	number = {3},
	pages = {884-909},
	year = {2019},
}

@article{affili,
	author = {E. Affili and E. Zuazua},
	title = {Controllability of diffusive Lotka–Volterra strongly competitive systems under boundary
constrained controls},
	journal = {preprint},
	volume = {},
	number = {},
	pages = {},
	year = {2026},
}

@article{bliman,
	author = {P.-A. Bliman and D. Cardona-Salgado and Y. Dumont and O. Vasilieva},
	journal = {Journal of Mathematical Sciences},
	number = {5},
	pages = {607-622},
	title = {Optimal control approach for implementation of sterile insect techniques},
	volume = {279},
	year = {2024},
}

\end{document}